\documentclass[12pt, reqno]{amsart}
\usepackage{amsmath, amstext, amsbsy, amssymb, amscd}

\newtheorem{theorem}{Theorem}[section]
\newtheorem{lemma}[theorem]{Lemma}
\newtheorem{proposition}[theorem]{Proposition}

\theoremstyle{definition}
\newtheorem{definition}[theorem]{Definition}
\newtheorem{example}[theorem]{Example}

\theoremstyle{remark}

\numberwithin{equation}{section}
\newcommand{\C}{ \mathbb C }

\newcommand{\End}{{\rm End}}
\newcommand{\fock}{{\mathbb H}_S}

\newcommand{\hil}[1]{{S^{[#1]}}}
\newcommand{\Hn}{H^*(\Sn)}

\newcommand{\Supp}{{\rm Supp}}

\newcommand{\vac}{|0\rangle}
\newcommand{\w}{\tilde}
\newcommand{\wa}{\tilde {\mathfrak a}}
\newcommand{\wt}{\tilde {\mathfrak t}}
\newcommand{\W}{\widetilde}
\newcommand{\Wq}{\widetilde Q}
\newcommand{\Wfock}{\widetilde {\mathbb H}_S}
\newcommand{\Sn}{S^{[n]}}

\newcommand{\Z}{ \mathbb Z }

{\vskip-\lastskip\medskip
  \noindent
  {\em #1.}\enspace
  }%
{\qed\par\medskip
  }

\begin{document}
\title[Incidence Hilbert schemes and Lie algebras]
      {Incidence Hilbert schemes and infinite dimensional Lie algebras}

\author[Wei-Ping Li]{Wei-Ping Li$^1$}
\address{Department of Mathematics, HKUST, Clear Water Bay, Kowloon,
Hong Kong} \email{mawpli@ust.hk}
\thanks{${}^1$Partially supported by the grant CERG601905}

\author[Zhenbo Qin]{Zhenbo Qin$^2$}
\address{Department of Mathematics, University of Missouri, Columbia,
MO 65211, USA} \email{zq@math.missouri.edu}
\thanks{${}^2$Partially supported by an NSF grant}

%\subjclass{Primary: 14C05; Secondary: 14F43, 17B65.}

\begin{abstract}
Let $S$ be a smooth projective surface. Using correspondences, we
construct an infinite dimensional Lie algebra that acts on the
direct sum
$$\Wfock=\bigoplus_{m=0}^{+\infty}H^*(S^{[m,m+1]})$$
of the cohomology groups of the incidence Hilbert schemes
$\hil{m,m+1}$. The algebra is related to an extension of an
infinite
%%change
dimensional Heisenberg algebra. The space $\Wfock$ is a highest
weight  representation of this algebra. Our result provides a
representation-theoretic interpretation of Cheah's generating
function of Betti numbers of the incidence Hilbert schemes. As a
consequence, an additive basis of $H^*(S^{[m,m+1]})$ is obtained.
\end{abstract}

\maketitle
\date{}
%\tableofcontents

%%
%%
%%
%%
%%
%%
\section{Introduction}
Let $S$ be a smooth projective surface. The Hilbert scheme $S^{[n]}$
of $n$ points on $S$ parametrizes all the $0$-dimensional closed
subschemes of $S$ with length $n$. It is a crepant resolution of the
$n$-fold symmetric product $S^{(n)}$. The generating function of the
Betti numbers of $S^{[n]}$ has an elegant closed formula due to
G\"ottsche \cite{Got1}:
\begin{eqnarray*}
\sum_{n=0}^{+ \infty}\left (\sum_{i = 0}^{4n} (-1)^i
   b_i\big ( \hil{n} \big ) z^i\right )
   q^n=\prod_{n=1}^{+\infty}
   \prod_{i = 0}^{4} \left ( \frac{1}{1-z^{2n-2+i}q^n}
   \right )^{(-1)^i b_i(S)},
\end{eqnarray*}
where $b_i(\cdot)$ denotes the $i$-th Betti number. There is a
representation theoretic interpretation discovered by Grojnowski
\cite{Gro} and Nakajima \cite{Na1} independently. It says that the
infinite dimensional vector space
\begin{eqnarray*}
\mathbb H_S=\bigoplus_{n=0}^{+ \infty}H^*(S^{[n]})
\end{eqnarray*}
is the Fock space of an infinite dimensional Heisenberg algebra. The
algebra is constructed geometrically via correspondences. This
result not only shows a beautiful symmetry for the Hilbert schemes
$S^{[n]}$, but also provides the right platform for the intensive
research on the cohomology ring of $S^{[n]}$ (e.g. see \cite{CG,
Got2, Lehn, LS1, LS2, LQW1, LQW2, LQW3, LQW4, LQW5, Mar1, Mar2, OP,
QW1, QW2, Vas}).

These Hilbert schemes $S^{[n]}$ play special roles in the recent
study of interplay between Donaldson-Thomas theory and Gromov-Witten
theory \cite{MNOP1, MNOP2, OP, KLQ, EQ}. When $S$ is embedded in a
smooth projective threefold $X$, the moduli spaces for the relative
Donaldson-Thomas theory of the pair $(X, S)$ admit natural morphisms
to the Hilbert schemes $S^{[n]}$. In some cases \cite{KLQ, EQ}, the
moduli spaces for the (ordinary) Donaldson-Thomas theory admit
morphisms to the Hilbert schemes $S^{[n]}$ as well. It is natural to
speculate whether certain interesting Lie algebra, perhaps somehow
related to the previously-mentioned Heisenberg algebra of the
Hilbert schemes, acts on the cohomologies of these moduli spaces of
the (ordinary and relative) Donaldson-Thomas theory.

We are making preliminary progress in this direction \cite{LQ2}. It
turns out that the incidence Hilbert schemes $S^{[n,n+1]}$ are part
of the ingredient in our work, where
\begin{eqnarray*}
\hil{n, n+1}=\{ (\xi,\xi')  \,|\, \xi \subset \xi' \} \subset
\hil{n}\times \hil{n+1}.
\end{eqnarray*}
These incidence Hilbert schemes $S^{[n,n+1]}$ are smooth. The
generating function of their Betti numbers has a closed formula due
to Cheah \cite{Ch2}:
\begin{eqnarray}  \label{Hnn1_bi}
& &\sum_{n=0}^{+\infty} \left (\sum_{i = 0}^{4(n+1)} (-1)^i
      b_i\big ( \hil{n, n+1} \big ) z^i\right )q^n    \nonumber  \\
&=&\left (\sum_{i=0}^4 (-1)^i b_i(S) z^i \right ) \cdot
   \frac{1}{1-z^2q} \cdot  \prod_{n=1}^{+\infty}
   \prod_{i = 0}^4 \left (\frac{1}{1-z^{2n-2+i}q^n}
   \right )^{(-1)^i b_i(S)}.
\end{eqnarray}
The infinite product in (\ref{Hnn1_bi}) hints that there should be a
Heisenberg algebra action, similar to the case of the Hilbert
schemes $S^{[n]}$, on the space
\begin{eqnarray*}
\Wfock = \bigoplus_{n=0}^{+\infty} H^* \big ( \hil{n, n+1} \big ).
\end{eqnarray*}
The two other factors in the formula (\ref{Hnn1_bi}) suggest that a
larger Lie algebra be needed for which the space $\Wfock$ is a
highest weight module.

In this paper, we show that such an infinite dimensional Lie algebra
exists, whose character formula agrees with (\ref{Hnn1_bi}), and
that the space $\Wfock$ is the highest weight representation of this
Lie algebra. In the following, we outline the basic ingredients in
the construction of this Lie algebra and its representation on
$\Wfock$.

 For $m \ge 0$ and $n > 0$, we define certain closed subset
$\Wq^{[m+n,m]}$ of
\begin{eqnarray*}
\hil{m+n, m+n+1} \times S \times \hil{m, m+1}.
\end{eqnarray*}
 For $\alpha \in H^*(S)$, we define
the operator $\wa_{-n}(\alpha) \in \End \big (\Wfock \big )$ by
\begin{eqnarray*}
\wa_{-n}(\alpha)(\W A) = \tilde{p}_{1*}([\Wq^{[m+n,m]}] \cdot
\tilde{\rho}^*\alpha \cdot \tilde{p}_2^*\W A)
\end{eqnarray*}
for $\W A \in H^*(\hil{m, m+1})$, where $\tilde{p}_1, \tilde{\rho},
\tilde{p}_2$ are the projections of
\begin{eqnarray*}
\hil{m+n, m+n+1} \times S \times \hil{m, m+1}
\end{eqnarray*}
to $\hil{m+n, m+n+1}, S, \hil{m, m+1}$ respectively. Define
$\wa_{n}(\alpha) \in \End \big (\Wfock \big )$ to be $(-1)^n$ times
the operator obtained from the definition of $\wa_{-n}(\alpha)$ by
switching the roles of $\tilde{p}_1$ and $\tilde{p}_2$. We define
$\wa_0(\alpha) = 0$ for every $\alpha \in H^*(S)$.

The first result is that the operators $\wa_{n}(\alpha), n \in \Z,
\alpha \in H^*(S)$ satisfy the following Heisenberg algebra
commutation relation:
\begin{eqnarray*}
[\wa_n(\alpha), \wa_k(\beta)] = -n \; \delta_{n,-k} \int_S(\alpha
\beta) \cdot {\rm Id}_{\Wfock}.
\end{eqnarray*}
This  explains the appearance of the infinite product in the formula
(\ref{Hnn1_bi}).

The two other factors in (\ref{Hnn1_bi}) correspond to two new
features of the space $\Wfock$. The first feature comes from an
$H^*(S)$-module structure of $\Wfock$. We observe that sending a
pair $(\xi,\xi') \in \hil{n, n+1}$ to the support of
$I_\xi/I_{\xi'}$ yields a morphism:
\begin{eqnarray*}
\rho_n: \hil{n, n+1} \to S.
\end{eqnarray*}
It follows that $H^*\big ( \hil{n, n+1} \big )$ and hence $\Wfock$
are $H^*(S)$-modules. It turns out that the operator
$\wa_{n}(\alpha): \Wfock \to \Wfock$ is an $H^*(S)$-module
homomorphism  for all $n$. This explains the appearance of the first
factor in the formula (\ref{Hnn1_bi}).

The second feature is the existence of a new operator $\wt \in \End
\big (\Wfock \big )$, called translation operator. It is constructed
via a correspondence of $\Wq_n$ which is a closed subset of
$S^{[n+1, n+2]} \times S^{[n,n+1]}$. It maps the component
$H^*(S^{[n,n+1]})$ to the component $H^*(S^{[n+1, n+2]})$, and is
commutative with the Heisenberg algebra, i.e.,
\begin{eqnarray*}
[\wt, \wa_{n}(\alpha)] = 0
\end{eqnarray*}
for all $n$ and $\alpha$. In addition, $\wt$ is an $H^*(S)$-module
homomorphism. The translation operator $\wt$ is responsible for the
second factor in (\ref{Hnn1_bi}).

Now we can construct the Lie algebra $\W{\mathfrak H}_S$ and state
our main theorem. Let ${\w {\mathfrak h}}_S$ be the Heisenberg
algebra generated by the operators $\wa_n(\alpha)$, $n \in \Z$,
$\alpha \in H^*(S)$ and the identity operator ${\rm Id}_{\Wfock}$.
Define a Lie algebra structure on
\begin{eqnarray*}
\W{\mathfrak H}_S = {\w {\mathfrak h}}_S \oplus H^*(S)\oplus \C\wt
%%change
\end{eqnarray*}
compatible with the Lie algebra structure on ${\w {\mathfrak h}}_S$
and the cup product on $H^*(S)$.
\begin{theorem}   \label{int_thm}
The space $\Wfock$ is a representation of the Lie algebra
$\W{\mathfrak H}_S $ with a highest weight vector being the vacuum
vector
\begin{eqnarray*}
\vac = 1_S \in H^0(S) = H^0(\hil{0,1})
\end{eqnarray*}
 where $1_S$ denotes the fundamental cohomology class of $S$.
 %%change
\end{theorem}

The organization of the paper is as follows. In \S
\ref{sect_basics}, we review the Heisenberg algebra action on the
Fock space $\mathbb H_S$ for the Hilbert schemes $S^{[n]}$ and some
basics of the incidence Hilbert schemes $S^{[n,n+1]}$. In \S
\ref{sect_action_inc}, we construct the Heisenberg operators
$\wa_n(\alpha)$ on $\Wfock$ and verify the Heisenberg commutation
relation. We also compare $\wa_n(\alpha)$ with the pullback of the
Heisenberg operator $\mathfrak a_n(\alpha)$ on $\mathbb H_S$ via the
morphism $f_n$ from $S^{[n,n+1]}$ to $S^{[n]}$. In \S
\ref{sect_trans}, we define the translation operator $\wt$, and show
that it commutes with the Heisenberg operators and is an
$H^*(S)$-module homomorphism.  We also compare $\wa_n(\alpha)$ with
the pullback of $\mathfrak a_n(\alpha)$ via the morphism $g_{n+1}$
from $S^{[n,n+1]}$ to $S^{[n+1]}$. It turns out that the
understanding of the pullback of the creation operators via
$g_{n+1}$ needs the translation operator. In \S \ref{sect_loop}, we
prove Theorem~\ref{int_thm}.

\medskip\noindent
{\bf Conventions.} Unless otherwise indicated, all the cohomology
groups in this paper are in $\C$-coefficients. For a continuous map
$p: Y_1 \to Y_2$ between two smooth compact manifolds and for
$\alpha_1 \in H^*(Y_1)$, we define $p_*(\alpha_1)$ to be ${\rm
PD}^{-1}p_{*}({\rm PD}(\alpha_1))$ where ${\rm PD}$ stands for the
Poincar\'e duality. The $\hbox{ }\widetilde{}\hbox{ }$ is used for
all the notations related to the incidence Hilbert schemes
$S^{[n,n+1]}$. When a product $X_1\times X_2\times\ldots X_n$ is
clear from the context, we use $\pi_{i_1 \cdots i_k}$ to denote the
projection from $X_1\times X_2\times\ldots X_n$ to the product of
the $i_1$-th, \ldots, $i_k$-th factors.

\bigskip\noindent
{\bf Acknowledgments.} The first author would like to thank the
Department of Mathematics at the University of Missouri for the
Miller Scholarship which made his visit there in February of 2006
possible and MSRI at Berkeley for its support. We thank Weiqiang
Wang for valuable suggestions.

\section{\bf Basics on Hilbert schemes of points on surfaces}
\label{sect_basics}

\subsection{\bf Hilbert schemes of points on surfaces}
\label{subsect_hil} $\,$ Let $S$ be a complex smooth projective
surface, and $\Sn$ be the Hilbert scheme of points in $S$. An
element in $\Sn$ is represented by a length-$n$ $0$-dimensional
closed subscheme $\xi$ of $S$. For $\xi \in \Sn$, let $I_{\xi}$ be
the corresponding sheaf of ideals. It is well known that $\Sn$ is
smooth. Sending an element in $\Sn$ to its support in the symmetric
product ${\rm Sym}^n(S)$, we obtain the Hilbert-Chow morphism
$\pi_n: \Sn \rightarrow {\rm Sym}^n(S)$, which is a resolution of
singularities. We have the universal codimension-$2$ subscheme:
\begin{eqnarray}  \label{cod2}
{\mathcal Z}_n=\{(\xi, s) \subset \Sn\times S \, | \, s\in {\rm
Supp}{(\xi)}\}\subset \Sn\times S.
\end{eqnarray}
Let $\Hn$ be the total cohomology of $\Sn$ with $\C$-coefficients.
Put
\begin{eqnarray}  \label{fock}
\fock = \bigoplus_{n=0}^{+\infty} \Hn.
\end{eqnarray}

For $m \ge 0$ and $n > 0$, let $Q^{[m,m]} = \emptyset$ and define
$Q^{[m+n,m]}$ to be the closed subset:
\begin{eqnarray*}
\{ (\xi, s, \eta) \in S^{[m+n]} \times S \times S^{[m]} \, | \, \xi
\supset \eta \text{ and } \mbox{Supp}(I_\eta/I_\xi) = \{ s \} \}.
\end{eqnarray*}
We recall Nakajima's definition of the Heisenberg operators
\cite{Na1, Na2}. Let $n > 0$. The linear operator $\mathfrak
a_{-n}(\alpha) \in \End(\fock)$ with $\alpha \in H^*(S)$ is defined
by
\begin{eqnarray}  \label{def_a}
\mathfrak a_{-n}(\alpha)(A) = p_{1*}([Q^{[m+n,m]}] \cdot
\rho^*\alpha \cdot p_2^*A)
\end{eqnarray}
for $A \in H^*(S^{[m]})$, where $p_1, \rho, p_2$ are the projections
of $S^{[m+n]} \times S \times S^{[m]}$ to $S^{[m+n]}, S, S^{[m]}$
respectively. Define $\mathfrak a_{n}(\alpha) \in \End(\fock)$ to be
$(-1)^n$ times the operator obtained from the definition of
$\mathfrak a_{-n}(\alpha)$ by switching the roles of $p_1$ and
$p_2$. We often refer to $\mathfrak a_{-n}(\alpha)$ (resp.
$\mathfrak a_n(\alpha)$) as the {\em creation} (resp. {\em
annihilation})~operator. We also set $\mathfrak a_0(\alpha) =0$. A
non-degenerate super-symmetric bilinear form $\langle -, - \rangle$
on $\mathbb H_S$ is induced from the standard one on $H^*(\Sn)$
defined by
\begin{eqnarray*}
\langle A, B \rangle  =\int_{\hil{n}} A B, \qquad A, B \in
H^*(\hil{n}).
\end{eqnarray*}
This allows us to define the {\it adjoint} $\mathfrak f^\dagger \in
\End(\mathbb H_S)$ for $\mathfrak f \in \End(\mathbb H_S)$. Then,
\begin{eqnarray}
&\mathfrak a_n(\alpha) = (-1)^n \cdot \mathfrak
a_{-n}(\alpha)^\dagger,&             \label{a_n_pm}  \\
&(\mathfrak f \circ \mathfrak g)^\dagger = \mathfrak g^\dagger \circ
\mathfrak f^\dagger.&     \label{fg}
\end{eqnarray}

It was proved in \cite{Na1, Na2} that the operators $\mathfrak
a_n(\alpha), n \in \Z, \alpha \in H^*(S)$ satisfy the following
Heisenberg algebra commutation relation:
\begin{eqnarray}  \label{eq:heis}
\displaystyle{[\mathfrak a_n(\alpha), \mathfrak a_k(\beta)] = -n \;
\delta_{n,-k} \int_S(\alpha \beta) \cdot {\rm Id}_{\fock}}.
\end{eqnarray}
Moreover, the space $\fock$ is an irreducible module over the
Heisenberg algebra generated by the operators $\mathfrak
a_n(\alpha)$ with a highest weight vector
\begin{eqnarray*}
\vac=1 \in H^0(S^{[0]}) \cong \C.
\end{eqnarray*}
It follows that $\fock$ is linearly spanned by all the {\it
Heisenberg monomial classes}:
\begin{eqnarray}  \label{Heis_mon}
\mathfrak a_{-n_1}(\alpha_1) \cdots \mathfrak a_{-n_k}(\alpha_k)\vac
\end{eqnarray}
where $k \ge 0, n_1, \ldots, n_k > 0$, and $\alpha_1, \ldots,
\alpha_k$ run over a linear basis of $H^*(S)$. We remark that the
Lie brackets in (\ref{eq:heis}) are understood in the super sense
according to the parity of the degrees of the cohomology classes
involved.

\subsection{\bf Incidence Hilbert schemes}
\label{subsect_inc} $\,$ Define the incidence variety:
\begin{eqnarray}  \label{Hnn1}
\hil{n, n+1}=\{ (\xi,\xi')  \,|\, \xi \subset \xi' \} \subset
\hil{n}\times \hil{n+1}.
\end{eqnarray}
It is well-known (see \cite{Ch1, ES}) that the incidence Hilbert
scheme $\hil{n, n+1}$ is irreducible, smooth and of dimension
$2(n+1)$, and that $\hil{n,n+1}$ is also the blowup of $S^{[n]}
\times S$ along the universal codimension-$2$ subscheme ${\mathcal
Z}_n$.

Besides the Hilbert scheme $\hil{n}$ for a surface $S$, the
incidence Hilbert scheme $\hil{n,n+1}$ for a surface $S$ is the only
class of (generalized or nested) Hilbert schemes of points on smooth
varieties of dimension bigger than one which are smooth for all $n$
(see \cite{Ch1}). Inspired by the results on the Hilbert schemes
$\hil{n}$, it is natural to consider the sum of the total cohomology
group of $S^{[n, n+1]}$ over all $n\ge 0$:
\begin{eqnarray}  \label{wfock}
\Wfock = \bigoplus_{n=0}^{+\infty} H^* \big ( \hil{n, n+1} \big ).
\end{eqnarray}

The space $\Wfock$ has richer structures than $\mathbb H_S$. For
instance, $\Wfock$ is an $H^*(S)$-module. The module structure is
induced by the morphism:
\begin{eqnarray}   \label{rho_n}
\rho_n: \hil{n, n+1} \to S
\end{eqnarray}
sending a pair $(\xi,\xi') \in \hil{n, n+1}$ to the support of
$I_\xi/I_{\xi'}$.

The following facts on punctual Hilbert schemes are useful for later
sections. Fix a point $s \in S$. For $m \ge 0$ and $n > 0$, we
define two closed subsets:
\begin{eqnarray}
   M_m(s)
&=&\{ \xi \in \hil{m} |\, \Supp(\xi) = \{s\}\},
               \label{mms}  \\
   M_{m, m+n}(s)
&=&\{ (\xi, \xi') |\, \xi \subset \xi' \}
   \subset M_m(s) \times M_{m+n}(s).  \label{mmns}
\end{eqnarray}
It is known that $M_{m, m+1}(s)$ and $M_{m+1}(s)$ are irreducible
with
\begin{eqnarray}  \label{dim_mms}
\dim M_{m, m+1}(s) = \dim M_{m+1}(s) = m.
\end{eqnarray}

\section{\bf Heisenberg algebra actions for incidence Hilbert schemes}
\label{sect_action_inc}

\subsection{\bf General remarks about correspondences.}
\label{subsect_cor} $\,$ Given two smooth projective varieties $X_1$
and $X_2$ and a closed subset $Z$ of $X_1\times X_2$, we can define
a map $[Z]_*$, called the correspondence of $Z$, from $H^*(X_2)$ to
$H^*(X_1)$ via
\begin{eqnarray*}
Z_*(A)=p_{1*}([Z]\cdot p_2^*(A))
\end{eqnarray*}
for $A\in H^*(X_2)$, where $p_i$ is the projection from $X_1\times
X_2$ to its $i$-th factor. Given another closed subset $Y\subset
X_2\times X_3$, the composition
\begin{eqnarray*}
[Z]_*\circ [Y]_*\colon H^*(X_3)\to H^*(X_1)
\end{eqnarray*}
is given by the correspondence of $(\pi_{13})_*(\pi_{12}^*[Z] \cap
\pi_{23}^*[Y])$. Note that we have used the notations established in
the Conventions.

In this paper, we will exhibit a collection of correspondences for
the incidence Hilbert schemes whose induced operators form a Lie
algebra. From the paragraph above, we see that to check the
commutation relations, we need to study
\begin{eqnarray*}
\pi_{13}\big (\pi_{12}^{-1}(Z)\cap \pi_{23}^{-1}(Y) \big ).
\end{eqnarray*}
This subset may not be irreducible, and the intersection
$\pi_{12}^{-1}(Z)\cap \pi_{23}^{-1}(Y)$ may not be transversal.
However, in most of the cases,
\begin{eqnarray*}
\pi_{12}^{-1}(Z) \cap \pi_{23}^{-1}(Y)
\end{eqnarray*}
has only one or two components with the expected dimensions, while
other components have smaller dimensions and thus will not
contribute to the induced map on cohomologies. In addition, either
the transversal property is satisfied, or the multiplicities can be
computed. Therefore, most of our work is to calculate the dimensions
of various stratum. Also, instead of doing the computations for all
the cases, we will only present the most exemplary ones in full
detail while lay out what the remaining cases are with details
skipped.

Even though a large percentage of our work is on the computation of
dimensions of various subsets, one should not be mislead to believe
that this is the key. It is free to construct the correspondence map
between the cohomologies of $X$ and $Y$ for any closed subset of
$X\times Y$. However, very few collections of correspondences
provide meaningful operators on cohomologies. To make things worse,
there is no guiding principle for choosing the ``right" ones in most
situations. In fact, the Heisenberg operators to be studied below
were found with hints from our study of another type of moduli
spaces for Donaldson-Thomas invariants. Once the ``right"
correspondences are found, the checking of the commutation relations
will be relatively straight forward, although it may be very
lengthy.

\subsection{\bf Definition of the Heisenberg operators}
\label{subsect_def_optr} $\,$ For $m \ge 0$ and $n > 0$, let
\begin{eqnarray*}
\Wq^{[m,m]} = \emptyset
\end{eqnarray*}
and define $\Wq^{[m+n,m]}$ to be the following closed subset of
$\hil{m+n, m+n+1} \times S \times \hil{m, m+1}$:
\begin{eqnarray*}
\Wq^{[m+n,m]}= \big \{ \big ( (\xi,\xi'), s, (\eta,\eta') \big )
|&&\xi
\supset \eta, \xi' \supset \eta', \quad \Supp(I_\eta/I_\xi)= \{ s \}, \\
&&\Supp(I_\xi/I_{\xi'}) = \Supp(I_\eta/I_{\eta'}) \big \}.
\end{eqnarray*}

The dimension of the subset $\Wq^{[m+n,m]}$ is given by the
following Lemma.

\begin{lemma}  \label{dim_Wq}
For $m \ge 0$ and $n > 0$, $\dim \Wq^{[m+n,m]} = 2m + n + 3$.
\end{lemma}
\begin{proof}
Take an element $\big ( (\xi,\xi'), s, (\eta,\eta') \big )$ in
$\Wq^{[m+n,m]}$. Let
\begin{eqnarray*}
\{t\}=\Supp(I_\xi/I_{\xi'})=\Supp(I_\eta/I_{\eta'})
\end{eqnarray*}
for some $t \in S$. First of all, assume that $s \ne t$. Then $\eta$
can be decomposed as
\begin{eqnarray*}
\eta = \eta_0+\eta_s + \eta_t
\end{eqnarray*}
where $\eta_s \in M_i(s)$ for some $i \ge 0$, $\eta_t \in M_j(t)$
for some $j \ge 0$, $\eta_0 \in \hil{m-i-j}$, and $s, t \not \in
\Supp(\eta_0)$. Then, $\eta', \xi$ and $\xi'$ can be written as
\begin{eqnarray*}
\eta' = \eta_0+\eta_s + \eta_t',   \quad \xi   = \eta_0+\xi_s +
\eta_t, \quad \xi'  = \eta_0+\xi_s + \eta_t',
\end{eqnarray*}
\begin{eqnarray*}
(\eta_t, \eta_t') \in M_{j,j+1}(t), \quad (\eta_s, \xi_s) \in
M_{i,i+n}(s).
\end{eqnarray*}
When $i = j = 0$, we conclude from (\ref{dim_mms}) that the number
of moduli of these triples $\big ( (\xi,\xi'), s, (\eta,\eta') \big
) \in \Wq^{[m+n,m]}$ is equal to
\begin{eqnarray}   \label{dim_Wq.1}
& &\# (\text{moduli of } \eta_0)
   + \# (\text{moduli of } \xi_s)
   + \# (\text{moduli of }  t)   \nonumber  \\
&=&2m + (n-1) + 4   \nonumber  \\
&=&2m + n + 3.
\end{eqnarray}

In general, when $i \ge 1$ or $j \ge 1$, we see from (\ref{dim_mms})
again that the number of moduli of these triples $\big ( (\xi,\xi'),
s, (\eta,\eta') \big ) \in \Wq^{[m+n,m]}$ is at most
\begin{eqnarray}   \label{dim_Wq.2}
& &\# (\text{moduli of } \eta_0) +
   \# (\text{moduli of } \eta_s) +
 \# (\text{moduli of } \eta_t \subset \eta_t')
   + \# (\text{moduli of } \xi_s)     \nonumber  \\
&=&2(m-i-j)+ {\rm max}(i-1, 0) + j  +(i+n-1)  + 4
   \nonumber  \\
&<&2m + n + 3.
\end{eqnarray}

Next, let $s = t$. This time we decompose the $0$-cycle $\eta$ into
\begin{eqnarray*}
\eta = \eta_0 + \eta_s
\end{eqnarray*}
where $\eta_s \in M_i(s)$ for some $i \ge 0$, $\eta_0 \in
\hil{m-i}$, and $s \not \in \Supp(\eta_0)$. Then,
\begin{eqnarray*}
\eta' = \eta_0 + \eta_s',   \quad \xi   &=& \eta_0 + \xi_s, \quad
\xi' = \eta_0 + \xi_s',
\end{eqnarray*}
\begin{eqnarray*}
(\eta_s, \eta_s') \in M_{i,i+1}(s),\quad&&(\eta_s, \xi_s) \in
M_{i,i+n}(s),\\
 (\xi_s, \xi_s') \in M_{i+n,i+n+1}(s),\quad
&&(\eta_s', \xi_s') \in M_{i+1,i+n+1}(s).
\end{eqnarray*}
So the number of moduli of these triples $\big ( (\xi,\xi'), s,
(\eta,\eta') \big ) \in \Wq^{[m+n,m]}$ is at most
\begin{eqnarray}   \label{dim_Wq.3}
& &\# (\text{moduli of } \eta_0) +
   \# (\text{moduli of } \eta_s \subset \eta_s') +
   \# (\text{moduli of } \xi_s \subset \xi_s') \nonumber  \\
&=&2(m-i) + i + (i+n) + 2                \nonumber  \\
&<&2m + n + 3.
\end{eqnarray}

Combining (\ref{dim_Wq.1}), (\ref{dim_Wq.2}) and (\ref{dim_Wq.3}),
we get $\dim \Wq^{[m+n,m]} = 2m + n + 3$.
\end{proof}

Now we introduce the operators on $\Wfock$ induced by the
correspondence of $\Wq^{[m+n,m]}$.
\begin{definition}  \label{def_wa}
(i) Let $\alpha \in H^*(S)$. Define $\wa_0(\alpha) = 0$ in $\End
\big (\Wfock \big )$.

(ii) Let $n > 0$ and $\alpha \in H^*(S)$. Define $\wa_{-n}(\alpha)
\in \End \big (\Wfock \big )$ by
\begin{eqnarray}  \label{def_wa.1}
\wa_{-n}(\alpha)(\W A) = \tilde{p}_{1*}([\Wq^{[m+n,m]}] \cdot
\tilde{\rho}^*\alpha \cdot \tilde{p}_2^*\W A)
\end{eqnarray}
for $\W A \in H^*(\hil{m, m+1})$, where $\tilde{p}_1, \tilde{\rho},
\tilde{p}_2$ are the projections of $\hil{m+n, m+n+1} \times S
\times \hil{m, m+1}$ to $\hil{m+n, m+n+1}, S, \hil{m, m+1}$
respectively.

(iii) Define $\wa_{n}(\alpha) \in \End \big (\Wfock \big )$ to be
$(-1)^n$ times the operator obtained from the definition of
$\wa_{-n}(\alpha)$ by switching the roles of $\tilde{p}_1$ and
$\tilde{p}_2$.
\end{definition}

A non-degenerate super-symmetric bilinear form $\langle -, -
\rangle$ on $\Wfock$ is induced from the standard one on
$H^*(\hil{n, n+1})$ defined by
\begin{eqnarray*}
\langle \W A, \W B \rangle  =\int_{\hil{n, n+1}} \W A \W B, \qquad
\W A, \W B \in H^*(\hil{n, n+1}).
\end{eqnarray*}
This allows us to define the {\it adjoint} $\w {\mathfrak f}^\dagger
\in \End \big (\Wfock \big )$ for $\w {\mathfrak f} \in \End \big
(\Wfock \big )$. Then,
\begin{eqnarray}   \label{wa_n_pm}
\wa_n(\alpha) = (-1)^n \cdot \wa_{-n}(\alpha)^\dagger.
\end{eqnarray}

Define $|\alpha| = \ell$ if $ \alpha \in H^\ell(S)$. Let $n \ne 0$.
By (\ref{def_wa.1}) and Lemma~\ref{dim_Wq}, we see that
$\wa_{-n}(\alpha)$ has bi-degree $(n, 2n-2 + |\alpha|)$, i.e.,
\begin{eqnarray*}
\wa_{-n}(\alpha): \, H^r \big (\hil{m, m+1} \big ) \to H^{r + 2n-2 +
|\alpha|} \big (\hil{m+n, m+n+1} \big ).
\end{eqnarray*}

\subsection{\bf $H^*(S)$-linearity}
\label{subsect_mod_homo} $\,$
 Recall from \S \ref{subsect_inc} that
$\Wfock$ is an $H^*(S)$-module.

\begin{lemma}  \label{wa_homo}
The map $\wa_{-n}(\alpha): \Wfock \to \Wfock$ is an $H^*(S)$-module
homomorphism.
\end{lemma}
\begin{proof}
We will only prove the Lemma for $n > 0$ since the proof for $n < 0$
is similar. Let $n > 0$, $\beta \in H^*(S)$, and $\W A \in H^* \big
(S^{[m, m+1]} \big )$. We need to show that
\begin{eqnarray}   \label{wa_homo.1}
\wa_{-n}(\alpha)(\rho_m^*\beta \cdot \W A) = (-1)^{|\alpha| |\beta|}
\,\, \rho_{m+n}^*\beta \cdot \wa_{-n}(\alpha)(\W A)
\end{eqnarray}
where the morphism $\rho_m^*$ is defined in (\ref{rho_n}). By
definition,
\begin{eqnarray*}
   \wa_{-n}(\alpha)(\rho_m^*\beta \cdot \W A)
&=&\tilde{p}_{1*} \big ([\Wq^{[m+n,m]}] \cdot
   \tilde{\rho}^*\alpha \cdot \tilde{p}_2^*(\rho_m^*\beta
   \cdot \W A) \big )     \\
&=&(-1)^{|\alpha| |\beta|} \,\, \tilde{p}_{1*} \big (
   [\Wq^{[m+n,m]}] \cdot \tilde{p}_2^*\rho_m^*\beta \cdot
   \tilde{\rho}^*\alpha \cdot \tilde{p}_2^*\W A \big )    \\
&=&(-1)^{|\alpha| |\beta|} \,\, \tilde{p}_{1*} \big (
   \tau_*((\rho_m \circ \tilde{p}_2 \circ \tau)^*\beta) \cdot
   \tilde{\rho}^*\alpha \cdot \tilde{p}_2^*\W A \big )
\end{eqnarray*}
where $\tau: \, \Wq^{[m+n,m]} \to \hil{m+n, m+n+1} \times S \times
\hil{m, m+1}$ is the inclusion map.

From the proof of Lemma~\ref{dim_Wq}, we see that for $\big (
(\xi,\xi'), s, (\eta,\eta') \big ) \in \Wq^{[m+n,m]}$,
\begin{eqnarray*}
\Supp(I_{\eta'}/I_{\xi'})= \Supp(I_{\eta}/I_{\xi})=\{ s \}.
\end{eqnarray*}
Therefore, $\rho_m \circ \tilde{p}_2 \circ \tau = \rho_{m+n} \circ
\tilde{p}_1 \circ \tau$. Hence, we obtain
\begin{eqnarray*}
   \wa_{-n}(\alpha)(\rho_m^*\beta \cdot \W A)
&=&(-1)^{|\alpha| |\beta|} \,\, \tilde{p}_{1*} \big (
   \tau_*((\rho_{m+n} \circ \tilde{p}_1 \circ \tau)^*\beta)
   \cdot \tilde{\rho}^*\alpha \cdot \tilde{p}_2^*\W A \big )\\
&=&(-1)^{|\alpha| |\beta|} \,\, \tilde{p}_{1*} \big (
   [\Wq^{[m+n,m]}] \cdot \tilde{p}_1^*\rho_{m+n}^*\beta
   \cdot \tilde{\rho}^*\alpha \cdot \tilde{p}_2^*\W A \big )\\
&=&(-1)^{|\alpha| |\beta|} \,\, \rho_{m+n}^*\beta \cdot
   \tilde{p}_{1*} \big ( [\Wq^{[m+n,m]}] \cdot
   \cdot \tilde{\rho}^*\alpha \cdot \tilde{p}_2^*\W A \big )
\end{eqnarray*}
by the projection formula. It follows immediately that
(\ref{wa_homo.1}) holds.
\end{proof}

\subsection{\bf Comparisons with Heisenberg operators on $\mathbb H_S$}
\label{subsect_com} $\,$
 There are two natural morphisms from
$\hil{m, m+1}$ to $\hil{m}$ and $\hil{m+1}$ respectively:
\begin{eqnarray*}
\begin{array}{ccc}
\hil{m, m+1}&\overset {g_{m+1}} \longrightarrow&\hil{m+1}\\
{\,\,\,} \downarrow {f_m}&&\\
\hil{m}.&&
\end{array}
\end{eqnarray*}

Since there are Heisenberg operators $\mathfrak a_{-n}(\alpha)$ on
$\mathbb H_S$, it is natural to ask whether the following diagram is
commutative or not:
\begin{eqnarray}\label{com_diag}
\CD H^* \big (S^{[m]} \big ) @>{\mathfrak a_{-n}(\alpha)}>>
    H^* \big (S^{[m+n]} \big )
\\ @VV{f_m^*}V @VV{f_{m+n}^*}V \\
H^* \big (\hil{m, m+1} \big ) @>{\wa_{-n}(\alpha)}>>
    H^* \big (\hil{m+n, m+n+1} \big ).
\endCD
\end{eqnarray}

\begin{lemma}  \label{wa_comm_f}
The diagram (\ref{com_diag}) is commutative .
\end{lemma}

This Lemma not only relates different Heisenberg operators on
different spaces, but also is used in the proof of
Proposition~\ref{prop_hei_inc} to determine a constant.

\begin{proof} Again, we will only prove the Lemma for $n > 0$.

Let $n > 0$ and $A \in H^* \big (S^{[m]} \big )$.  By (\ref{def_a})
and (\ref{def_wa.1}), we obtain
\begin{eqnarray}
   f_{m+n}^*\mathfrak a_{-n}(\alpha)(A)
&=&f_{m+n}^*p_{1*}([Q^{[m+n,m]}] \cdot
       \rho^*\alpha \cdot p_2^*A),    \label{wa_comm_f.1} \\
   {\wa_{-n}(\alpha)}f_m^*(A)
&=&\tilde{p}_{1*}([\Wq^{[m+n,m]}] \cdot
   \tilde{\rho}^*\alpha \cdot \tilde{p}_2^*f_m^*(A)).
                                      \label{wa_comm_f.2}
\end{eqnarray}
Let $g = {\rm Id}_{\hil{m+n, m+n+1} \times S} \times f_m$ and $h =
f_{m+n} \times {\rm Id}_{S \times \hil{m}}$:
\begin{eqnarray*}
&g:\; &\hil{m+n, m+n+1} \times S \times \hil{m, m+1} \to
         \hil{m+n, m+n+1} \times S \times \hil{m},  \\
&h:\; &\hil{m+n, m+n+1} \times S \times \hil{m} \to
         \hil{m+n} \times S \times \hil{m},
\end{eqnarray*}
and let ${\w p}_1': \hil{m+n, m+n+1} \times S \times \hil{m} \to
\hil{m+n, m+n+1}$ be the first projection. By (\ref{wa_comm_f.1})
and base change, we conclude that
\begin{eqnarray}
   f_{m+n}^*\mathfrak a_{-n}(\alpha)(A)
&=&({\w p}_1')_* h^*([Q^{[m+n,m]}] \cdot
   \rho^*\alpha \cdot p_2^*A)                 \nonumber \\
&=&({\w p}_1')_* \big (h^*[Q^{[m+n,m]}] \cdot
   (\rho \circ h)^*\alpha \cdot (p_2 \circ h)^*A \big ).
                                   \label{wa_comm_f.3}
\end{eqnarray}
 Similarly, using (\ref{wa_comm_f.2}), $\w
p_1=\w p_1^{\prime}\circ g$, and the projection formula, we obtain
\begin{eqnarray}
   {\wa_{-n}(\alpha)}f_m^*(A)
&=&({\w p}_1')_* g_* \left ([\Wq^{[m+n,m]}] \cdot
   g^* \big ( (\rho \circ h)^*\alpha \cdot
   (p_2 \circ h)^*A \big ) \right )         \nonumber \\
&=&({\w p}_1')_* \big (g_*[\Wq^{[m+n,m]}] \cdot
   (\rho \circ h)^*\alpha \cdot (p_2 \circ h)^*A \big ).
                                   \label{wa_comm_f.4}
\end{eqnarray}

Next, we compare the two cycles $h^*[Q^{[m+n,m]}]$ and
$g_*[\Wq^{[m+n,m]}]$. Note that both $h^{-1}(Q^{[m+n,m]})$ and
$g(\Wq^{[m+n,m]})$ are equal to the subset
\begin{eqnarray*}
W :=\left \{ \big ( (\xi,\xi'), s, \eta \big ) |\, \xi \supset \eta
\text{ and } \Supp(I_\eta/I_\xi) = \{ s \} \right \}
\end{eqnarray*}
of $\hil{m+n, m+n+1} \times S \times \hil{m}$. A typical element
$(\xi, s, \eta) \in Q^{[m+n,m]}$ is of the form:
\begin{eqnarray*}
(s_1 + \ldots + s_{m} + \xi_s, \,\, s, \,\, s_1 + \ldots + s_{m})
\end{eqnarray*}
where $s_1, \ldots, s_{m}, s \in S$ are distinct, and $\xi_s \in
M_n(s)$. Hence,
\begin{eqnarray}       \label{wa_comm_f.5}
h^*[Q^{[m+n,m]}] = [W].
\end{eqnarray}
On the other hand, a typical element $\big ( (\xi,\xi'), s, (\eta,
\eta') \big ) \in \Wq^{[m+n,m]}$ is of the form:
\begin{eqnarray*}
\big (\left (\xi_0+ \xi_s, \, \xi_0 +s_{m+1} + \xi_s \right ), \,\,
s, \,\, \left (\xi_0, \, \xi_0 +s_{m+1}\right ) \big )
\end{eqnarray*}
where $s_1, \ldots, s_{m+1}, s \in S$ are distinct, $\xi_0=
s_1+\ldots+s_m $, and $\xi_s \in M_n(s)$. Hence,
\begin{eqnarray}       \label{wa_comm_f.6}
g_*[\Wq^{[m+n,m]}] = [W].
\end{eqnarray}

The combinations of  (\ref{wa_comm_f.3}), (\ref{wa_comm_f.4}),
(\ref{wa_comm_f.5}) and (\ref{wa_comm_f.6}) prove the Lemma.
\end{proof}

We remark that in \S \ref{subsect_ant_comp}, another comparison of
the Heisenberg operators $\mathfrak a_{-n}(\alpha)$ and
$\wa_{-n}(\alpha)$ by using the pullbacks $g_{m+1}^*$ will be
presented.

\subsection{\bf Heisenberg commutation relation}
\label{subsect_Heisenberg_rela} $\,$ The proof of the following
Proposition is a bit long. It is divided into many cases. If one
understands one case, one should be able to follow the proof of
remaining cases easily.

\begin{proposition}  \label{prop_hei_inc}
The operators $\wa_{n}(\alpha), n \in \Z, \alpha \in H^*(S)$ satisfy
the following Heisenberg algebra commutation relation:
\begin{eqnarray}  \label{prop_hei_inc.1}
[\wa_n(\alpha), \wa_k(\beta)] = -n \; \delta_{n,-k} \int_S(\alpha
\beta) \cdot {\rm Id}_{\Wfock}.
\end{eqnarray}
\end{proposition}
\begin{proof}
In view of (\ref{wa_n_pm}), formula (\ref{prop_hei_inc.1}) is
equivalent to the formulas:
\begin{eqnarray}
   [\wa_{-n}(\alpha), \wa_{-k}(\beta)]
&=&0,                        \label{prop_hei_inc.2}   \\
   {[\wa_{-n}(\alpha), \wa_{k}(\beta)]}
&=&n \; \delta_{n,k} \int_S(\alpha \beta) \cdot
   {\rm Id}_{\Wfock}         \label{prop_hei_inc.3}
\end{eqnarray}
where $n, k > 0$. We will prove these two formulas separately.

\medskip
\noindent{\it Proof of (\ref{prop_hei_inc.2})}:
\label{subsect_pf_hei_inc.2} $\,$
\par

Let $n, k > 0$. Then the operator $\wa_{-n}(\alpha) \wa_{-k}(\beta)$
is induced by the class:
\begin{eqnarray}  \label{class_w}
w = \pi_{1245*} \left ( \pi^*_{123} \big [ \Wq^{[m+n+k, m+k]} \big ]
\cdot \pi^*_{345} \big [ \Wq^{[m+k,m]} \big ] \right ) \in
A_{2m+k+n+4}(W')
\end{eqnarray}
where $A_*(\cdot)$ denotes the Chow group, and $W' = \pi_{1245}(W)$
with
\begin{eqnarray}  \label{corrs1}
W = \pi_{123}^{-1} \big ( \Wq^{[m+n+k, m+k]} \big ) \cap
\pi_{345}^{-1} \big ( \Wq^{[m+k,m]} \big ).
\end{eqnarray}

Note that $W$ is a closed subset of the ambient space
\begin{eqnarray*}
 S^{[m+n+k,
m+n+k+1]}\times S \times S^{[m+k, m+k+1]} \times S\times S^{[m,
m+1]}.
\end{eqnarray*}
In the sequel, the ambient spaces in similar situations will not be
explicitly presented since they can be written out easily from the
context.

 From Lemma~\ref{dim_Wq}, we know that the expected dimension
of the intersection $W$ should be $2m+n+k+4$. The subset $W$ may
have many irreducible components. Those that are mapped by
$\pi_{1245}$ to subsets of dimension less than $2m+n+k+4$ will not
contribute to the cycle $w$. The aim of the computations below is to
pick out those components with dimension no less than the expected
dimension.

Any element $((\xi, \xi^{\prime}), s, (\eta, \eta^{\prime}), t,
(\zeta, \zeta^{\prime}))$ in $W$ must satisfy the conditions:
\begin{eqnarray*}
&\zeta\subset \zeta^{\prime},\quad \eta\subset\eta^{\prime}, \quad
\xi\subset\xi^{\prime}, \quad \zeta\subset \eta \subset \xi,\quad
\zeta^{\prime}\subset
\eta^{\prime}\subset \xi^{\prime},&   \\
&\Supp(I_{\eta}/I_{\xi})=\{ s \},
\quad \Supp(I_{\zeta}/I_{\eta})=\{ t \},&  \\
&\Supp(I_{\zeta}/I_{\zeta^{\prime}})
=\Supp(I_{\eta}/I_{\eta^{\prime}}) =\Supp(I_{\xi}/I_{\xi^{\prime}})
=\{ p \}&
\end{eqnarray*}
for some $p \in S$. In the following, we consider four different
cases separately. We use $W_i$ to denote the subset of $W$
consisting of all the points satisfying Case $i$, put
\begin{eqnarray*}
\mathfrak w_i = \dim \, \pi_{1245}(W_i)
\end{eqnarray*}
(similar notations such as $U_i, \mathfrak u_i, V_i, \mathfrak v_i$
will be used throughout the rest of the paper).

\medskip\noindent{\bf Case 1}: $s, t, p$ are distinct.
We have the following decompositions:
\begin{eqnarray*}
\zeta=\zeta_0+\zeta_s+\zeta_t+\zeta_p,&&
\zeta^{\prime}=\zeta_0+\zeta_s+\zeta_t+\zeta_p^{\prime},\\
\eta=\zeta_0+\zeta_s+\eta_t+\zeta_p,&&
\eta^{\prime}=\zeta_0+\zeta_s+\eta_t+\zeta_p^{\prime},\\
\xi=\zeta_0+\xi_s+\eta_t+\zeta_p,&&
\xi^{\prime}=\zeta_0+\xi_s+\eta_t+\zeta_p^{\prime},
\end{eqnarray*}
where $\Supp(\zeta_0)\cap\{s, t, p\}=\emptyset$, $\ell(\zeta_s)=i$,
$\ell(\zeta_t)=j$, $\ell(\zeta_p)=\ell$, and $\zeta_s, \zeta_t,
\ldots$ are supported at $\{s\}, \{t\}, \ldots$ respectively. Then,
\begin{eqnarray*}
   \mathfrak w_1
&=&\#(\text{moduli of }\zeta_0)+
     \#(\text{moduli of }\zeta_s)+
     \#(\text{moduli of }\zeta_t)+\\
& &+ \, \#(\text{moduli of }\eta_t) +
     \#(\text{moduli of }\xi_s) +
     \#(\text{moduli of }\zeta_p\subset \zeta_p^{\prime}) \\
&=&2(m-i-j-\ell)+ {\rm max}(i-1, 0)+2+ {\rm max}(j-1, 0) +2+ \\
& &+ \,(j+k-1)+(i+n-1)+\ell+2  \\
&=&(2m+k+n+4) - i + {\rm max}(i-1, 0) -j
     + {\rm max}(j-1, 0) -\ell.
\end{eqnarray*}

\noindent{\it Case 1.1}: If $i=j=\ell=0$, then we have $\mathfrak
w_1 = (2m+k+n+4)$.

\medskip\noindent
{\it Case 1.2}: If one of the integers $i, j$ and $\ell$ is
positive, then $\mathfrak w_1 < (2m+k+n+4)$.

\medskip
The three remaining cases are listed by:
\begin{enumerate}
\item[] {\bf Case 2:} $s=t\neq p$;

\item[] {\bf Case 3:} $s=p\neq t$ (by symmetry, this also covers
the case $t = p \ne s$);

\item[] {\bf Case 4:} $s=p=t$.
\end{enumerate}
For these three cases, we skip the arguments which are in the same
style as in Case 1. They all have dimension less than the expected
dimension $(2m+n+k+4)$. Therefore, only Case 1.1 has contribution to
the cohomological operation. In this subcase, it is not difficult to
show that the intersection (\ref{corrs1}) along $W_1$ is
transversal. Moreover, $\pi_{1245}(W_1)$ consists of all the points
of the form:
\begin{eqnarray*}
\big ((\zeta_0+\xi_s+\eta_t, \zeta_0+\xi_s+\eta_t+  p ), s, t,
(\zeta_0, \zeta_0 + p ) \big )
\end{eqnarray*}
in $S^{[m+n+k, m+n+k+1]}\times S \times S\times S^{[m, m+1]}$. So
the contribution of this subcase to the operator $\wa_{-n}(\alpha)
\wa_{-k}(\beta)$ coincides with the corresponding contribution for
the operator $(-1)^{|\alpha| |\beta|} \cdot \wa_{-k}(\beta)
\wa_{-n}(\alpha)$. In other words, we obtain the identity:
\begin{eqnarray*}
\wa_{-n}(\alpha) \wa_{-k}(\beta) = (-1)^{|\alpha| |\beta|} \cdot
\wa_{-k}(\beta) \wa_{-n}(\alpha).
\end{eqnarray*}
This completes the proof of the commutation relation
(\ref{prop_hei_inc.2}).

\medskip
\noindent{\it Proof of (\ref{prop_hei_inc.3})}:
\label{subsect_pf_hei_inc.3} $\,$
\par

Let $n, k > 0$. Then the operator $\wa_{-n}(\alpha) \wa_{k}(\beta)$
is induced by the class:
\begin{eqnarray}  \label{class_u}
u = \pi_{1245*} \left ( \pi^*_{123} \big [ \Wq^{[m+n, m]} \big ]
\cdot \pi^*_{345} \big [ \tau_{m, m+k} (\Wq^{[m+k, m]}) \big ]
\right ) \in A_{2m+k+n+4}(U')
\end{eqnarray}
where $\tau_{m, m+k}: \, S^{[m+k, m+k+1]} \times S \times S^{[m,
m+1]} \to S^{[m, m+1]} \times S \times S^{[m+k, m+k+1]}$ is the
isomorphism switching $S^{[m, m+1]}$ and $S^{[m+k, m+k+1]}$, and $U'
= \pi_{1245}(U)$ with
\begin{eqnarray}  \label{corrs2}
U = \pi_{123}^{-1} \big ( \Wq^{[m+n, m]} \big ) \cap \pi_{345}^{-1}
\big ( \tau_{m, m+k} (\Wq^{[m+k, m]}) \big ).
\end{eqnarray}

Any element $((\xi, \xi^{\prime}), s, (\eta, \eta^{\prime}), t,
(\zeta, \zeta^{\prime}))$ in $U$  must satisfy the conditions
\begin{eqnarray*}
&\zeta\subset \zeta^{\prime},\quad \eta\subset\eta^{\prime},\quad
\xi\subset\xi^{\prime}, \quad \xi\supset \eta\subset \zeta, \quad
\xi^{\prime}
\supset \eta^{\prime}\subset \zeta^{\prime},&     \\
&\Supp(I_{\eta}/I_{\xi})=\{s\},
\quad \Supp(I_{\eta}/I_{\zeta})=\{t\},&  \\
&\Supp(I_{\zeta}/I_{\zeta^{\prime}})
=\Supp(I_{\eta}/I_{\eta^{\prime}}) =\Supp(I_{\xi}/I_{\xi^{\prime}})
=\{ p \}&
\end{eqnarray*}
for some $p \in S$. In the following, we consider four different
cases separately.

\medskip\noindent{\bf Case 1}: $s, t, p$ are distinct.
We have the following decompositions:
\begin{eqnarray*}
\eta=\eta_0+\eta_s+\eta_t+\eta_p,&&
\eta^{\prime}=\eta_0+\eta_s+\eta_t+\eta_p^{\prime},\\
\xi=\eta_0+\xi_s+\eta_t+\eta_p,&&
\xi^{\prime}=\eta_0+\xi_s+\eta_t+\eta_p^{\prime},\\
\zeta=\eta_0+\eta_s+\zeta_t+\eta_p,&&
\zeta^{\prime}=\eta_0+\eta_s+\zeta_t+\eta_p^{\prime}
\end{eqnarray*}
where $\Supp(\eta_0)\cap\{s, t, p\}=\emptyset$, $\ell(\eta_s)=i$,
$\ell(\eta_t)=j$, and $\ell(\eta_p)=\ell$. Then,
\begin{eqnarray*}
\mathfrak u_1 &=&\#(\text{moduli of }\eta_0)+\#(\text{moduli of
}\xi_s)
   +\#(\text{moduli of }\zeta_t)+ \\
& &+ \, \#(\text{moduli of }\eta_s) +
   \#(\text{moduli of }\eta_t)+
   \#(\text{moduli of }\eta_p\subset \eta_p^{\prime})  \\
&=&2(m-i-j-\ell)+(i+n-1)+(j+k-1)+ \\
& &+\, {\rm max}(i-1, 0)+{\rm max}(j-1, 0)+\ell+6\\
&=&(2m+k+n+4) - i + {\rm max}(i-1, 0) -j
     + {\rm max}(j-1, 0) -\ell.
\end{eqnarray*}

\noindent{\it Case 1.1}: If $i=j=\ell=0$, then we have $\mathfrak
u_1 = (2m+k+n+4)$.

\medskip\noindent
{\it Case 1.2}: If one of the integers $i, j$ and $\ell$ is
positive, then $\mathfrak u_1 < (2m+k+n+4)$.

\medskip
The three remaining cases are listed by:
\begin{enumerate}
\item[] {\bf Case 2:} $s=t\neq p$;

\item[] {\bf Case 3:} $s=p\neq t$ (by symmetry, this also covers
the case $t = p \ne s$);

\item[] {\bf Case 4:} $s=p=t$.
\end{enumerate}
In these three cases, all the dimensions are smaller than the
expected dimension $(2m+n+k+4)$. So only Case 1.1 contributes to the
class $u$ in (\ref{class_u}).

\medskip
Next we consider the operator $\wa_{k}(\beta)\wa_{-n}(\alpha)$. This
is the only case where there are two components for $V$ below with
the expected dimension. One of the components will cancel out with
the one from $\wa_{-n}(\alpha)\wa_{k}(\beta)$, and the other is a
non-transversal intersection which carries a multiplicity . Using
Lemma~\ref{wa_comm_f} which compares $\wa_{n}(\alpha)$ with the
Heisenberg operators $\mathfrak a_{n}(\alpha)$ on $\mathbb H_S$, we
determine the multiplicity.

More precisely, the operator $\wa_{k}(\beta)\wa_{-n}(\alpha)$ is
induced by the class:
\begin{eqnarray}  \label{class_v}
   v
= \pi_{1245*} \left ( \pi^*_{123} \big [
   \tau_{m+n, m+k+n} (\Wq^{[m+k+n, m+n]}) \big ] \cdot
   \pi^*_{345} \big [ \Wq^{[m+k+n, m+k]} \big ] \right )
\end{eqnarray}
in $A_{2m+k+n+4}(V')$, where $V' = \pi_{1245}(V)$ and $V$ is given
by
\begin{eqnarray}  \label{corrs3}
V = \pi_{123}^{-1} \big ( \tau_{m+n, m+k+n} (\Wq^{[m+k+n, m+n]})
\big ) \cap \pi_{345}^{-1} \big ( \Wq^{[m+k+n, m+k]} \big ).
\end{eqnarray}
Any element $((\xi, \xi^{\prime}), s, (\eta, \eta^{\prime}), t,
(\zeta, \zeta^{\prime}))$ in (\ref{corrs3}) must satisfy the
conditions:
\begin{eqnarray*}
&\zeta\subset \zeta^{\prime},\quad \eta\subset\eta^{\prime}, \quad
\xi\subset\xi^{\prime}, \quad \xi\subset \eta\supset \zeta,\quad
\xi^{\prime}\subset \eta^{\prime}\supset \zeta^{\prime},&  \\
&\Supp(I_{\xi}/I_{\eta})=\{s\},
\quad \Supp(I_{\zeta}/I_{\eta})=\{ t\},&    \\
&\Supp(I_{\zeta}/I_{\zeta^{\prime}})
=\Supp(I_{\eta}/I_{\eta^{\prime}}) =\Supp(I_{\xi}/I_{\xi^{\prime}})
=\{ p \}&
\end{eqnarray*}
for some $p \in S$. In the following, we consider four different
cases separately.

\medskip\noindent{\bf Case 1${}^\prime$}: $s, t, p$ are distinct.
We have the following decompositions:
\begin{eqnarray*}
\zeta=\zeta_0+\zeta_s+\zeta_t+\zeta_p,&&
\zeta^{\prime}=\zeta_0+\zeta_s+\zeta_t+\zeta_p^{\prime},\\
\eta=\zeta_0+\zeta_s+\eta_t+\zeta_p,&&
\eta^{\prime}=\zeta_0+\zeta_s+\eta_t+\zeta_p^{\prime},\\
\xi=\zeta_0+\xi_s+\eta_t+\zeta_p,&&
\xi^{\prime}=\zeta_0+\xi_s+\eta_t+\zeta_p^{\prime},
\end{eqnarray*}
where $\Supp(\zeta_0)\cap\{s, t, p\}=\emptyset$, $\ell(\xi_s)=i$,
$\ell(\zeta_t)=j$, and $\ell(\zeta_p)=\ell$. Then,
\begin{eqnarray*}
\mathfrak v_1 &=&\#(\text{moduli of }\zeta_0)+\#(\text{moduli of
}\zeta_s)+
   \#(\text{moduli of }\zeta_t)+  \\
& &+ \, \#(\text{moduli of }\eta_t)+\#(\text{moduli of }\xi_s) +
   \#(\text{moduli of }\zeta_p\subset \zeta_p^{\prime})   \\
&=&2(m-i-j-\ell)+(i+k-1)+{\rm max}(j-1, 0) + \\
& &+\, (j+n-1) + {\rm max}(i-1, 0)+\ell+6\\
&=&(2m+k+n+4) - i + {\rm max}(i-1, 0) -j
     + {\rm max}(j-1, 0) -\ell.
\end{eqnarray*}

\noindent{\it Case 1${}^\prime$.1}: If $i=j=\ell=0$, then we have
$\mathfrak v_1 = (2m+k+n+4)$.

\medskip\noindent
{\it Case 1${}^\prime$.2}: If one of the integers $i, j$ and $\ell$
is positive, then $\mathfrak v_1 < (2m+k+n+4)$.

\medskip\noindent{\bf Case 2${}^\prime$}: $s=t\neq p$.
We have the following decompositions:
\begin{eqnarray*}
\zeta=\zeta_0+\zeta_s+\zeta_p,&&
\zeta^{\prime}=\zeta_0+\zeta_s+\zeta_p^{\prime},\\
\eta=\zeta_0+\eta_s+\zeta_p,&&
\eta^{\prime}=\zeta_0+\eta_s+\zeta_p^{\prime},\\
\xi=\zeta_0+\xi_s+\zeta_p,&&
\xi^{\prime}=\zeta_0+\xi_s+\zeta_p^{\prime},
\end{eqnarray*}
where $\Supp(\zeta_0)\cap\{s, p\}=\emptyset$, $\ell(\xi_s)=i$,
$\ell(\zeta_s)=j$, and $\ell(\zeta_p)=\ell$. Note that
\begin{eqnarray}   \label{ki=nj}
k+i=\ell(\eta_s)=n+j.
\end{eqnarray}

\noindent{\it Case 2${}^\prime$.1}: $i=j=\ell=0$. By (\ref{ki=nj}),
we have $k=n$. Thus,
\begin{eqnarray*}
\mathfrak v_2 = 2(m+k)+4=2m+k+n+4.
\end{eqnarray*}

\noindent{\it Case 2${}^\prime$.2}: If one of the integers $i, j$
and $\ell$ is positive, then $ \mathfrak v_2 <  (2m+k+n+4)$.

\medskip
The two remaining cases are listed by:
\begin{enumerate}
\item[] {\bf Case 3${}^\prime$:} $s=p\neq t$ (by symmetry, this also
covers the case $t = p \ne s$);

\item[] {\bf Case 4${}^\prime$:} $s=p=t$.
\end{enumerate}
In these cases, the dimensions are smaller than $(2m+k+n+4)$.
Therefore these two cases have no contributions to the class $v$ in
(\ref{class_v}).

\smallskip
Finally, note that the contribution of Case 1.1 to the class $u$ in
(\ref{class_u}) and the contribution of Case 1${}^\prime$.1 to the
class $v$ in (\ref{class_v}) cancel out in the commutation:
\begin{eqnarray}
[\wa_{-n}(\alpha), \wa_{k}(\beta)]. \label{prop_hei_inc.3'}
\end{eqnarray}
So $[\wa_{-n}(\alpha), \wa_{k}(\beta)] = 0$ when $n \ne k$. This
proves (\ref{prop_hei_inc.3}) when $n \ne k$.

When $n = k$, Case 2${}^\prime$.1 also contributes to the operator
$\wa_{k}(\beta)\wa_{-n}(\alpha)$, and hence to
(\ref{prop_hei_inc.3'}). Note from Case 2${}^\prime$.1 that
$\pi_{1245}(V_2)$ consists of all the points of the form:
\begin{eqnarray*}
\big ((\zeta_0, \zeta_0+  p ), s, s, (\zeta_0, \zeta_0 +  p ) \big )
\end{eqnarray*}
in $S^{[m+n, m+n+1]}\times S \times S\times S^{[m+n, m+n+1]}$, where
$s \ne p$ and $p \not \in \Supp(\zeta_0)$. Thus,
\begin{eqnarray*}
[\wa_{-n}(\alpha), \wa_{n}(\beta)]|_{H^*(S^{[m+n, m+n+1]})} \,\, =
\,\, c \cdot \int_S(\alpha \beta) \cdot {\rm Id}_{H^*(S^{[m+n,
m+n+1]})}
\end{eqnarray*}
for some constant $c$. Now we conclude from Lemma~\ref{wa_comm_f}
and (\ref{eq:heis}) that $c = n$. This completes the proof of the
commutation relation (\ref{prop_hei_inc.3}) when $n = k$.
\end{proof}

\section{\bf A translation operator for incidence Hilbert schemes}
\label{sect_trans}

In this section, we will introduce a new operator $\wt$ on $\Wfock$,
called the translation operator. The operator $\wt$ is constructed
via a correspondence and has many nice properties. Indeed, we will
show that it is an $H^*(S)$-module homomorphism, commutes with the
Heisenberg operators, and has a left inverse. These properties imply
that it is responsible for the second factor in (\ref{Hnn1_bi}).

From Lemma~\ref{wa_comm_f}, we may tend to infer that the Heisenberg
operators $\wa_n(\alpha)$ are no different from the pullback of the
more well-known Heisenberg operators $\mathfrak a_n(\alpha)$. This
is almost the case when we compare these two types of operators via
the map
\begin{eqnarray*}
f_m\colon \hil{m,m+1}\to \hil{m}.
\end{eqnarray*}
However, we will see that when we compare them via the map
\begin{eqnarray*}
g_{m+1}\colon \hil{m,m+1}\to \hil{m+1},
\end{eqnarray*}
they are far from the same. The difference somehow is measured by
the new translation operator. Once we are presented with the fact
that the space $\Wfock$ is a highest weight  module of the algebra
generated by the operators $\wa_n(\alpha)$, $\beta\in H^*(S)$ and
$\wt$, we realize that the naive choice of the pullback of
Heisenberg algebras on Hilbert schemes either by the map $f_m$ or
$g_{m+1}$ won't provide the right algebra. In this regard, the new
algebra we are going to construct is subtler and richer than the
Heisenberg algebra on the Hilbert scheme $\hil{m}$. The fundamental
difference between these two algebras is the translation operator
$\wt$.
\subsection{\bf Definition of the translation operator}
\label{subsect_def_trans} $\,$ Let $\Wq_m$ be the closed subset:
\begin{eqnarray*}
   \Wq_m
&=&\{\big ( (\xi^{\prime}, \xi^{\prime \prime}), (\xi,
   \xi^{\prime}) \big ) |\Supp(I_\xi/I_{\xi'}) =
   \Supp(I_{\xi'}/I_{\xi''})\}    \\
&\subset&S^{[m+1, m+2]} \times S^{[m,m+1]}.
\end{eqnarray*}
As in Lemma~\ref{dim_Wq}, we can verify that $\dim \Wq_m = 2m + 3$.

\begin{definition}  \label{def_oper_wt}
Define the linear operator $\wt \in \End \big (\Wfock \big )$ by
\begin{eqnarray}  \label{def_oper_wt.1}
\wt(\W A) = \tilde{p}_{1*}([\Wq_m] \cdot \tilde{p}_2^*\W A)
\end{eqnarray}
for $\W A \in H^*(\hil{m, m+1})$, where $\tilde{p}_1, \tilde{p}_2$
are the two projections of $S^{[m+1, m+2]} \times S^{[m,m+1]}$.
\end{definition}

The bi-degrees of $\wt$ and its adjoint $\wt^\dagger$ are $(1, 2)$
and $(-1, -2)$ respectively.

\subsection{\bf $H^*(S)$-linearity and the left inverse}
Recall that $\Wfock$ is an $H^*(S)$-module.

\begin{lemma}  \label{wt_homo}
{\rm (i)} The operator $(-\wt^\dagger)$ is the left inverse of
$\wt$;

{\rm (ii)} The maps $\wt, \wt^\dagger: \Wfock \to \Wfock$ are
$H^*(S)$-module homomorphisms.
\end{lemma}
\begin{proof}
We skip the proof of (ii) since it is similar to the proof of
Lemma~\ref{wa_homo}. In the following, we prove (i). Note that the
operator $\wt^\dagger \, \wt$ is induced by the class:
\begin{eqnarray}  \label{class_wt_adj}
w = \pi_{13*} \left ( \pi^*_{12} \big [ \tau(\Wq_m) \big ] \cdot
\pi^*_{23} \big [ \Wq_m \big ] \right ) \in A_{2m+2}(W')
\end{eqnarray}
where $\tau: \, S^{[m+1, m+2]} \times S^{[m, m+1]} \to S^{[m, m+1]}
\times S^{[m+1, m+2]}$ is the isomorphism switching the two factors,
$W' = \pi_{13}(W)$, and $W$ is given by
\begin{eqnarray}  \label{corrs_t}
W = \pi_{12}^{-1} \big ( \tau(\Wq_m) \big ) \cap \pi_{23}^{-1} \big
( \Wq_m \big ).
\end{eqnarray}
The expected dimension of the intersection $W$ is $2m+2$.

Any element $\big ( (\eta, \xi^{\prime}), (\xi^{\prime}, \xi^{\prime
\prime}), (\xi, \xi^{\prime}) \big )$ in $W$ must satisfy the
conditions:
\begin{eqnarray*}
&\eta, \xi \subset\xi^{\prime} \subset
\xi^{\prime\prime}, &   \\
&\Supp(I_{\eta}/I_{\xi^{\prime}})
=\Supp(I_{\xi^{\prime}}/I_{\xi^{\prime\prime}})
=\Supp(I_{\xi}/I_{\xi^{\prime}}) =\{ p \}&
\end{eqnarray*}
for some $p \in S$. We have the following decompositions:
\begin{eqnarray*}
\xi=\xi_0 + \xi_p,&&
\xi^{\prime}=\xi_0+\xi_p^{\prime},\\
\xi^{\prime\prime}=\xi_0+\xi_p^{\prime\prime},&& \eta=\xi_0+\eta_p,
\end{eqnarray*}
where $p \not \in \Supp(\xi_0)$ and $\ell(\xi_p)=i$. Note that
$\ell(\eta_p)=i$. Fix the integer $i$.

If  $i > 0$, then the projection of the subset of $W$, consisting of
all the elements
\begin{eqnarray*}
\big ( (\eta, \xi^{\prime}), (\xi^{\prime}, \xi^{\prime \prime}),
(\xi, \xi^{\prime}) \big ),
\end{eqnarray*}
to $S^{[m, m+1]} \times S^{[m, m+1]}$ has dimension at most equal
to:
\begin{eqnarray*}
& &\#(\text{moduli of }\xi_0)+
     \#(\text{moduli of }\xi_p \subset \xi_p^{\prime})+
     \#(\text{moduli of }\eta_p)     \\
&=&2(m-i)+ i +(i-1) + 2  \\
&<&(2m+2).
\end{eqnarray*}
So the case $i > 0$ does not contribute to the class $w$ in
(\ref{class_wt_adj}).

If $i = 0$, then $\xi_p^{\prime} = p$, $\xi_p^{\prime\prime} \in
M_2(p)$, and $\xi_p = \eta_p = \emptyset$. So $\eta = \xi = \xi_0$.
The projection of this part of $W$ to $S^{[m, m+1]} \times S^{[m,
m+1]}$ consists of elements of the form:
\begin{eqnarray*}
\left \{ \big ( (\xi_0, \xi_0 + p), (\xi_0, \xi_0 + p) \big ) \right
\}.
\end{eqnarray*}
It follows that $w = c [\Delta]$ for some constant $c$, where
$\Delta$ denotes the diagonal in $S^{[m, m+1]} \times S^{[m, m+1]}$.
Hence, we get $\wt^\dagger \, \wt = c \, {\rm Id}_{H^*(S^{[m,
m+1]})}$.

To determine $c$, note that we can split off the factor $\xi_0$ from
our consideration, i.e., we can simply consider the case $m = 0$.
Then we have the morphism:
\begin{eqnarray*}
\pi_{13}: S \times \hil{1, 2} \times S \to S \times S.
\end{eqnarray*}
Fix a point $p \in S$. Since $[\Delta] \cdot [\{p\} \times S] = 1$,
we see from (\ref{class_wt_adj}) that
\begin{eqnarray*}
   c
&=&w \cdot [\{p\} \times S]       \\
&=&\pi_{13*} \left ( \pi^*_{12} \big [ \tau(\Wq_0) \big ]
   \cdot \pi^*_{23} \big [ \Wq_0 \big ] \right )
   \cdot [\{p\} \times S]       \\
&=&\pi^*_{12} \big [ \tau(\Wq_0) \big ]
   \cdot \pi^*_{23} \big [ \Wq_0 \big ]
   \cdot [\{p\} \times \hil{1, 2} \times S]    \\
&=&[\{p\} \times U_p \times S] \cdot [\{p\} \times \Wq_0]
\end{eqnarray*}
where $U_p = \big \{ (p, \xi_p^{\prime \prime})|\, \xi_p^{\prime
\prime} \in M_2(p) \big \}$. It follows that
\begin{eqnarray}  \label{wt_homo.1}
   c
= [U_p \times S] \cdot [\Wq_0] = [U_p] \cdot \phi_{1*}[\Wq_0] =
[U_p] \cdot [\phi_{1}(\Wq_0)]
\end{eqnarray}
where $\phi_1: \hil{1, 2} \times S \to \hil{1, 2}$ is the
projection. We have
\begin{eqnarray*}
\phi_{1}(\Wq_0) = \{ (s, \xi_s^{\prime \prime})|\, s \in S \text{
and } \xi_s^{\prime \prime} \in M_2(s) \big \}.
\end{eqnarray*}
Recall the natural morphism $g_2: \hil{1,2} \to \hil{2}$. It is
known from \cite{ES} that
\begin{eqnarray*}
[\phi_{1}(\Wq_0)] = \frac{1}{2} \, g_2^*[M_2(S)]
\end{eqnarray*}
where $M_2(S) = \cup_{s \in S} M_2(s)$. Hence, we see from
(\ref{wt_homo.1}) that
\begin{eqnarray*}
   c
&=&[U_p] \cdot \frac{1}{2} \, g_2^*[M_2(S)]
   = \frac{1}{2} \, g_{2*}[U_p] \cdot [M_2(S)]  \\
&=&\frac{1}{2} \, [g_2(U_p)] \cdot [M_2(S)]
   = \frac{1}{2} \, [M_2(p)] \cdot [M_2(S)]
   = -1
\end{eqnarray*}
where we have used the fact that $[M_2(p)] \cdot [M_2(S)] = -2$ from
\cite{ES}.
\end{proof}

\subsection{\bf Commutativity with the Heisenberg operators}
$\,$  The translation operator $\wt$ may look similar to the
creation operators at the first glimpse. However, we now show that
it differs from the creation operators in the essential way in that
it commutes with all the annihilation operators.

\begin{proposition}  \label{comm_wtwa} The translation operator
$\wt$ and its adjoint $\wt^\dagger$ commute with the Heisenberg
operators $\wa_{-n}(\alpha)$ for all $n$ and $\alpha$, i.e.,
\begin{eqnarray*}
[\wt, \wa_{-n}(\alpha)] =[\wt^\dagger, \wa_{-n}(\alpha)] = 0.
\end{eqnarray*}
\end{proposition}
\begin{proof}
Let $n > 0$. Then Proposition~\ref{comm_wtwa} is decomposed into two
parts:
\begin{eqnarray}
{[\wt, \wa_{-n}(\alpha)]}&=& 0, \label{comm_wtwa.1}\\
{[\wt, \wa_{n}(\alpha)]} &=& 0. \label{comm_wtwa.2}
\end{eqnarray}
We prove them separately, and will compare the proof with that of
Proposition~\ref{prop_hei_inc}.

\medskip\noindent
{\it Proof of (\ref{comm_wtwa.1})}:

Let $n > 0$. This part is similar to the proof of
(\ref{prop_hei_inc.2}).

The operator $\wa_{-n}(\alpha) \wt$ is induced by the class:
\begin{eqnarray}  \label{class_wt}
w = \pi_{124*} \left ( \pi^*_{123} \big [ \Wq^{[m+n+1, m+1]} \big ]
\cdot \pi^*_{34} \big [ \Wq_m \big ] \right ) \in A_{2m+n+4}(W')
\end{eqnarray}
where $W' = \pi_{124}(W)$ and the subset $W$ is defined by
\begin{eqnarray}  \label{corrs4}
W = \pi_{123}^{-1} \big ( \Wq^{[m+n+1, m+1]} \big ) \cap
\pi_{34}^{-1} \big ( \Wq_m \big ).
\end{eqnarray}
The expect dimension of the intersection $W$ is $2m+n+4$.

Any element $\big ( (\eta, \eta^{\prime}), s, (\xi^{\prime},
\xi^{\prime \prime}), (\xi, \xi^{\prime}) \big )$ in (\ref{corrs4})
must satisfy the conditions
\begin{eqnarray*}
&\eta\subset\eta^{\prime}, \quad \xi \subset\xi^{\prime} \subset
\xi^{\prime\prime}, \quad \xi^{\prime} \subset \eta, \quad
\xi^{\prime\prime} \subset \eta^{\prime},&   \\
&\Supp(I_{\xi^{\prime}}/I_{\eta})=\{ s \},&  \\
&\Supp(I_{\eta}/I_{\eta^{\prime}})
=\Supp(I_{\xi^{\prime}}/I_{\xi^{\prime\prime}})
=\Supp(I_{\xi}/I_{\xi^{\prime}}) =\{ p \}&
\end{eqnarray*}
for some $p \in S$. In the following, we consider two cases
separately.

\medskip\noindent{\bf Case 1}: $s \ne p$.
We have the following decompositions:
\begin{eqnarray*}
\xi=\xi_0+\xi_s+\xi_p,&& \xi^{\prime}=\xi_0+\xi_s+\xi_p^{\prime},
\quad \xi^{\prime\prime}=\xi_0+\xi_s+\xi_p^{\prime\prime},\\
\eta=\xi_0+\eta_s+\xi_p^{\prime},&&
\eta^{\prime}=\xi_0+\eta_s+\xi_p^{\prime\prime},
\end{eqnarray*}
where $\Supp(\xi_0)\cap\{s, p\}=\emptyset$, $\ell(\xi_s)=i$, and
$\ell(\xi_p)=j$.

\medskip\noindent{\it Case 1.1}: If $i=j=0$, then we have
$\mathfrak w_1 = (2m+n+4)$.

\medskip\noindent
{\it Case 1.2}: If $i> 0$ or $j > 0$, then we obtain $\mathfrak w_1
< (2m+n+4)$.

\medskip\noindent
{\bf Case 2}: $s = p$. We have $\mathfrak w_2 <(2m+n+4)$.
\medskip

In summary, only Case 1.1 contributes to the class $w$. In this
subcase, the intersection (\ref{corrs4}) along $W_1$ is transversal,
and $\pi_{1245}(W_1)$ consists of all the elements:
\begin{eqnarray*}
\big ((\xi_0 +\eta_s+  p , \xi_0 +\eta_s+ \xi_p^{\prime\prime}), s,
(\xi_0, \xi_0 +  p ) \big )
\end{eqnarray*}
in $S^{[m+n+1, m+n+2]}\times S \times S^{[m, m+1]}$. So the
contribution of this subcase to the operator $\wa_{-n}(\alpha) \,
\wt$ coincides with the corresponding contribution for the operator
$\wt \, \wa_{-n}(\alpha)$. This completes the proof of the
commutation relation (\ref{comm_wtwa.1}).

\medskip\noindent
{\it Proof of (\ref{comm_wtwa.2})}:

This part looks similar to the proof of (\ref{prop_hei_inc.3}).
However, it is fundamentally different. Recall that the subset $V$
in (\ref{corrs3}) whose correspondence defines
$\wa_{n}(\alpha)\wa_{-n}(\beta)$ has one more component with the
expected dimension than the subset for the operator
$\wa_{-n}(\alpha)\wa_{n}(\beta)$. It doesn't occur here. If there
were two components, they would come from the operator
$\wa_{n}(\alpha) \wt$. So pay a special attention to Case 2 below
and compare it with Case 2${}^{\prime}$.1 in the proof of
Proposition~\ref{prop_hei_inc}. In order to illustrate this subtle
difference, we present the full content of the proof.

Let $n > 0$. Then the operator $\wa_{n}(\alpha) \wt$ is induced by
the class:
\begin{eqnarray}  \label{class_ut}
u = \pi_{124*} \left ( \pi^*_{123} \big [ \tau_{m+1,
m+n+1}(\Wq^{[m+n+1, m+1]}) \big ] \cdot \pi^*_{34}\big [ \Wq_{m+n}
\big ] \right ) \in A_{2m+n+4}(U')
\end{eqnarray}
where $U' = \pi_{124}(U)$ and the subset $U$ is defined by
\begin{eqnarray}  \label{corrs5}
U = \pi_{123}^{-1} \big ( \tau_{m+1, m+n+1} (\Wq^{[m+n+1, m+1]})
\big ) \cap \pi_{34}^{-1} \big ( \Wq_{m+n} \big ).
\end{eqnarray}

Any element $\big ( (\eta, \eta^{\prime}), s, (\xi^{\prime},
\xi^{\prime \prime}), (\xi, \xi^{\prime}) \big )$ in (\ref{corrs5})
must satisfy the conditions
\begin{eqnarray*}
&\eta\subset\eta^{\prime}, \quad \xi \subset\xi^{\prime} \subset
\xi^{\prime\prime}, \quad \eta \subset \xi^{\prime}, \quad
\eta^{\prime} \subset \xi^{\prime\prime},&   \\
&\Supp(I_{\eta}/I_{\xi^{\prime}})=\{ s \},&  \\
&\Supp(I_{\eta}/I_{\eta^{\prime}})
=\Supp(I_{\xi^{\prime}}/I_{\xi^{\prime\prime}})
=\Supp(I_{\xi}/I_{\xi^{\prime}}) =\{ p \}&
\end{eqnarray*}
for some $p \in S$. In the following, we consider two cases
separately.

\medskip\noindent{\bf Case 1}: $s \ne p$.
We have the following decompositions:
\begin{eqnarray*}
\xi=\xi_0+\xi_s+\xi_p,&& \xi^{\prime}=\xi_0+\xi_s+\xi_p^{\prime},
\quad \xi^{\prime\prime}=\xi_0+\xi_s+\xi_p^{\prime\prime},\\
\eta=\xi_0+\eta_s+\xi_p^{\prime},&&
\eta^{\prime}=\xi_0+\eta_s+\xi_p^{\prime\prime},
\end{eqnarray*}
where $\Supp(\xi_0)\cap\{s, p\}=\emptyset$, $\ell(\eta_s)=i$, and
$\ell(\xi_p)=j$.

\medskip\noindent{\it Case 1.1}: If $i=j=0$, then we have
$\mathfrak u_1 = (2m+n+4)$.

\medskip\noindent
{\it Case 1.2}: If $i> 0$ or $j > 0$, then we obtain
\begin{eqnarray*}
     \mathfrak u_1
&\le&\#(\text{moduli of }\xi_0)+
     \#(\text{moduli of }\xi_s)+
     \#(\text{moduli of }\eta_s) + \\
&   &+ \, \#(\text{moduli of }\xi_p)+
     \#(\text{moduli of }\xi_p^{\prime}
     \subset \xi_p^{\prime\prime}) \\
&=  &2(m-i-j)+(i+n-1)+ {\rm max}(i-1, 0) + \\
&   &+ {\rm max}(j-1, 0)+ (j + 1)+4  \\
&<  &(2m+n+4).
\end{eqnarray*}

\noindent{\bf Case 2}: $s = p$. We have the following
decompositions:
\begin{eqnarray*}
\xi=\xi_0+\xi_s,&& \xi^{\prime}=\xi_0+\xi_s^{\prime},
\quad\xi^{\prime\prime}=\xi_0+\xi_s^{\prime\prime},\\
\eta=\xi_0+\eta_s,&&
 \eta^{\prime}=\xi_0+\eta_s^{\prime},
\end{eqnarray*}
where $s \not \in \Supp(\xi_0)$ and $\ell(\eta_s)=i$. In this case,
we see that
\begin{eqnarray*}
     \mathfrak u_2
&\le&\#(\text{moduli of }\xi_0)+
     \#(\text{moduli of }\xi_s \subset \xi_s^{\prime})+
     \#(\text{moduli of }\eta_s \subset \eta_s^{\prime}) \\
&\le  &2(m+1-i)+ (i+n-1)+ i  + 2     \\
&<  &(2m+n+4).
\end{eqnarray*}

In summary, only Case 1.1 contributes to the class $u$. The
contribution of this subcase to the operator $\wa_{n}(\alpha) \,
\wt$ coincides with that of the operator $\wt \, \wa_{n}(\alpha)$.
This completes the proof of the commutation relation
(\ref{comm_wtwa.2}).
\end{proof}

\subsection{\bf Another comparison with Heisenberg operators on
$\hil{m+1}$} \label{subsect_ant_comp} $\,$ This subsection is not
related to the Theorem in \S \ref{sect_loop}.  The main theme here
is to compare the Heisenberg operators $\wa_n(\alpha)$ with the
pull-back of the Heisenberg operators $\mathfrak a_n(\alpha)$ via
the morphism $g_{m+1}\colon \hil{m, m+1}\to \hil{m+1}$.

The comparison of annihilation operators is similar to
Lemma~\ref{wa_comm_f}. The proof of the following lemma is omitted
since it is similar to the proof of Lemma~\ref{wa_comm_f}.

\begin{lemma}  \label{wa_comm_g}
Let $n > 0$ and $\alpha \in H^*(S)$. Then, we have a commutative
diagram:
$$
\CD H^* \big (S^{[m+1]} \big ) @<{\mathfrak a_{n}(\alpha)}<<
    H^* \big (S^{[m+n+1]} \big )
\\ @VV{g_{m+1}^*}V @VV{g_{m+n+1}^*}V \\
H^* \big (\hil{m, m+1} \big ) @<{\wa_{n}(\alpha)}<<
    H^* \big (\hil{m+n, m+n+1} \big ).
\endCD
$$
\end{lemma}

Lemma~\ref{wa_comm_g} will not hold for the creation operators.
However, Proposition~\ref{wa_comm_g-n} below provides an explicit
formula relating the creation operators. In order to prove
Proposition~\ref{wa_comm_g-n}, we begin with a technical lemma.

Let $m \ge 0$. Define $\Wq_{m, 0}$ to be the diagonal of $S^{[m,
m+1]} \times S^{[m, m+1]}$. For $n \ge 1$, define $\Wq_{m, n}$ to be
the closed subset of $S^{[m+n, m+n+1]} \times S^{[m, m+1]}$:
\begin{eqnarray*}
\Wq_{m, n}=\big \{\big ( (\xi, \xi^{\prime}), (\eta, \eta^{\prime})
\big )|&&\xi^{\prime}\supset \xi\supset\eta^{\prime} \supset \eta,\\
&&\Supp(I_{\xi}/I_{\xi^{\prime}}) =\Supp(I_{\eta^{\prime}}/I_{\xi})
=\Supp(I_{\eta}/I_{\eta^{\prime}})\big \}.
\end{eqnarray*}
Note that $\Wq_{m, n}$ has dimension $(2m+n+2)$, and contains
exactly one irreducible component of dimension $(2m+n+2)$ whose
generic elements are of the form:
\begin{eqnarray}  \label{gen_pt_wq}
\big((\eta+\xi_s, \eta+\xi_s^{\prime}), (\eta, \eta+ s) \big ),
\qquad s \not \in \Supp(\eta).
\end{eqnarray}

\begin{lemma}    \label{wt^n}
The restriction of $\wt^n$ to $H^*(\hil{m, m+1})$ is given by the
cycle $[\Wq_{m, n}]$.
\end{lemma}
\begin{proof}
Use induction on $n$. Note that $\Wq_{m, 1} = \Wq_m$ by definition.
So the Lemma is trivially true when $n = 0, 1$. Next, assume that $n
\ge 2$ and the restriction of $\wt^{n-1}$ to $H^*(\hil{m, m+1})$ is
given by the cycle $[\Wq_{m, n-1}]$. Then we see that the
restriction of $\wt^{n} = \wt \, \wt^{n-1}$ to $H^*(\hil{m, m+1})$
is given by the cycle
\begin{eqnarray}    \label{wt^n.1}
     w
=\pi_{13*} \left ( \pi_{12}^*[\Wq_{m+n-1}]
     \cdot \pi_{23}^*[\Wq_{m, n-1}]\right )
\in A_{2m+n+2}(W')
\end{eqnarray}
where  $W' = \pi_{13}(W)$ and the subset $W$ is given by
\begin{eqnarray}    \label{wt^n.2}
W = \pi_{12}^{-1}(\Wq_{m+n-1}) \cap \pi_{23}^{-1}(\Wq_{m, n-1}).
\end{eqnarray}

Any element $\big((\xi, \xi^{\prime}), (\theta, \theta^{\prime}),
(\eta, \eta^{\prime}) \big )$ in (\ref{wt^n.2}) must satisfy:
\begin{eqnarray*}
&\xi^{\prime}\supset \xi=\theta^{\prime}\supset
 \theta\supset\eta^{\prime}\supset \eta,&\\
&\Supp(I_{\xi}/I_{\xi^{\prime}})
 =\Supp(I_{\theta}/I_{\theta^{\prime}})
 =\Supp(I_{\eta^{\prime}}/I_{\theta})
 =\Supp(I_{\eta}/I_{\eta^{\prime}}) = \{s \}&
\end{eqnarray*}
for some point $s \in S$. We have the following decompositions:
\begin{eqnarray*}
\eta=\eta_0+\eta_s,&& \eta^{\prime}=\eta_0+\eta_s^{\prime},\quad
\theta=\eta_0+\theta_s,\\
\theta^{\prime}=\xi=\eta_0+\xi_s  ,&&
\xi^{\prime}=\eta_0+\xi_s^{\prime},
\end{eqnarray*}
where $s \not \in \Supp(\eta_0)$ and $\eta_s \subset \eta_s^\prime
\subset \theta_s \subset \xi_s \subset \xi_s^{\prime}$. When $\ell
:= \ell(\eta_s) \ge 1$, the relation $\eta_s^\prime \subset \xi_s$
imposes a nontrivial condition on the points
\begin{eqnarray*}
\big ( (\eta_s, \eta_s^\prime), (\xi_s,\xi_s^{\prime})\big ) \in
M_{\ell, \ell+1}(s) \times M_{\ell+n, \ell+n+1}(s).
\end{eqnarray*}
Hence the subset consisting of the images $\big((\xi, \xi^{\prime}),
(\eta, \eta^{\prime}) \big )$ of those points
\begin{eqnarray*}
\big((\xi, \xi^{\prime}), (\theta, \theta^{\prime}), (\eta,
\eta^{\prime}) \big )
\end{eqnarray*}
in (\ref{wt^n.2}) with $\ell = \ell(\eta_s) \ge 1$ has dimension at
most
\begin{eqnarray*}
2(m - \ell) + 2 + [\ell+(n+\ell) - 1] = (2m+n+1)
\end{eqnarray*}
which is less than the expected dimension $(2m+n+2)$. So the case
$\ell \ge 1$ does not contribute to the cycle $w$ in (\ref{wt^n.1}).
When $\ell = 0$, we have
\begin{eqnarray*}
\eta=\eta_0,&& \eta^{\prime}=\eta_0+ s,\quad
\theta=\eta_0+\theta_s,\\
\theta^{\prime}=\xi=\eta_0+\xi_s  ,&&
\xi^{\prime}=\eta_0+\xi_s^{\prime}
\end{eqnarray*}
where $s \not \in \Supp(\eta_0)$ and $\theta_s \subset \xi_s \subset
\xi_s^{\prime}$. The images $\big((\eta_0+\xi_s,
\eta_0+\xi_s^{\prime}), (\eta_0, \eta_0+ s) \big )$ form a subset of
dimension $(2m+n+2)$ in $\hil{m+n,m+n+1} \times \hil{m,m+1}$. From
the description (\ref{gen_pt_wq}) of the generic points in $\Wq_{m,
n}$, we conclude that
\begin{eqnarray}     \label{wcWq}
w = c \cdot [\Wq_{m, n}]
\end{eqnarray}
where $c$ is the intersection multiplicity of (\ref{wt^n.2}) along
the generic points.

To determine $c$, choose a primitive integral class $\W A \in
H^*(\hil{m, m+1})$, i.e.,
\begin{eqnarray*}
\W A \in H^*(\hil{m, m+1}; \Z)/{\rm Tor} \subset H^*(\hil{m, m+1}).
\end{eqnarray*}
By (\ref{wcWq}), $\wt^n(\W A) = c \W B$ for some integral class $\W
B$. By Lemma~\ref{wt_homo}~(i), we see that
\begin{eqnarray*}
\W A = (-\wt^\dagger)^n \wt^n(\W A) = c (-\wt^\dagger)^n(\W B).
\end{eqnarray*}
Since $\W A$ is primitive and $(\wt^\dagger)^n(\W B)$ is an integral
class, we must have $c = 1$.
\end{proof}

\begin{proposition}  \label{wa_comm_g-n}
Let $n > 0$, $\alpha \in H^*(S)$, and $A \in H^* \big (S^{[m+1]}
\big )$. Then,

{\rm (i)} $g_{n}^* \mathfrak a_{-n}(\alpha)\vac = n \cdot
\wt^{n-1}\rho_{0}^*(\alpha)$;

{\rm (ii)} $g_{m+n+1}^* \mathfrak a_{-n}(\alpha)(A) =
\wa_{-n}(\alpha)(g_{m+1}^*A) + n \cdot \wt^{n-1} \big (
\rho_{m+1}^*(\alpha) \cdot f_{m+1}^*(A) \big )$.
\end{proposition}
\begin{proof}
(i) It suffices to prove the formula for $\alpha = 1_S$. The
cohomology class $\mathfrak a_{-n}(1_S)\vac$ is represented by the
subscheme
\begin{eqnarray*}
M_n(S) := \{ \xi^{\prime} \in \hil{n} |\, \Supp(\xi^{\prime}) =
\{s\} \,\, \text{\rm for some $s \in S$}\}.
\end{eqnarray*}
By the Lemma~2.4 in \cite{ES}, $1/n \cdot g_{n}^*\mathfrak
a_{-n}(1_S)\vac$ is represented by the subscheme
\begin{eqnarray*}
M_{n-1, n}(S) := \{ (\xi, \xi^{\prime}) \in \hil{n-1, n} |\,
\Supp(\xi^{\prime}) = \{s\} \,\, \text{\rm for some $s \in S$}\}.
\end{eqnarray*}
By Lemma~\ref{wt^n}, $\wt^{n-1}(1_{\hil{0, 1}})$ is also represented
by $M_{n-1, n}(S)$. Therefore,
\begin{eqnarray*}
g_{n}^* \mathfrak a_{-n}(1_S)\vac = n \cdot
\wt^{n-1}\rho_{0}^*(1_S).
\end{eqnarray*}

(ii) By the definitions (\ref{def_a}) and (\ref{def_wa.1}), we
obtain
\begin{eqnarray}
   g_{m+n+1}^*\mathfrak a_{-n}(\alpha)(A)
&=&g_{m+n+1}^*p_{1*}([Q^{[m+n+1,m+1]}] \cdot
       \rho^*\alpha \cdot p_2^*A),    \label{wa_comm_g-n.1} \\
   {\wa_{-n}(\alpha)}g_{m+1}^*(A)
&=&\tilde{p}_{1*}([\Wq^{[m+n,m]}] \cdot
   \tilde{\rho}^*\alpha \cdot \tilde{p}_2^*g_{m+1}^*(A)).
                                      \label{wa_comm_g-n.2}
\end{eqnarray}
Let $g = {\rm Id}_{\hil{m+n, m+n+1} \times S} \times g_{m+1}$ and $h
= g_{m+n+1} \times {\rm Id}_{S \times \hil{m+1}}$:
\begin{eqnarray*}
&g:\; &\hil{m+n, m+n+1} \times S \times \hil{m, m+1} \to
         \hil{m+n, m+n+1} \times S \times \hil{m+1},  \\
&h:\; &\hil{m+n, m+n+1} \times S \times \hil{m+1} \to
         \hil{m+n+1} \times S \times \hil{m+1}.
\end{eqnarray*}
Let ${\w p}_1': \hil{m+n, m+n+1} \times S \times \hil{m+1} \to
\hil{m+n, m+n+1}$ be the first projection. By (\ref{wa_comm_g-n.1})
and base change, we conclude that
\begin{eqnarray}
   g_{m+n+1}^*\mathfrak a_{-n}(\alpha)(A)
&=&({\w p}_1')_* h^*([Q^{[m+n+1,m+1]}] \cdot
   \rho^*\alpha \cdot p_2^*A)                 \nonumber \\
&=&({\w p}_1')_* \big (h^*[Q^{[m+n+1,m+1]}] \cdot
   (\rho \circ h)^*\alpha \cdot (p_2 \circ h)^*A \big ).
                                   \label{wa_comm_g-n.3}
\end{eqnarray}
Similarly, using (\ref{wa_comm_g-n.2}), $\w p_{1*}={\w p}_1'\circ g$
and the projection formula, we obtain
\begin{eqnarray}
   {\wa_{-n}(\alpha)}g_{m+1}^*(A)
&=&({\w p}_1')_* g_* \left ([\Wq^{[m+n,m]}] \cdot
   g^* \big ( (\rho \circ h)^*\alpha \cdot
   (p_2 \circ h)^*A \big ) \right )         \nonumber \\
&=&({\w p}_1')_* \big (g_*[\Wq^{[m+n,m]}] \cdot
   (\rho \circ h)^*\alpha \cdot (p_2 \circ h)^*A \big ).
                                   \label{wa_comm_g-n.4}
\end{eqnarray}

Next, we compare $h^*[Q^{[m+n+1,m+1]}]$ and $g_*[\Wq^{[m+n,m]}]$.
Let $U_1$ be the closure of the subset of $\hil{m+n, m+n+1} \times S
\times \hil{m+1}$ consisting of all the elements of the form:
\begin{eqnarray*}
\left (\left (\sum_{i=1}^m s_i + \xi'_s, \, \sum_{i=1}^{m+1} s_i +
\xi'_s \right ), \,\, s, \,\, \sum_{i=1}^{m+1} s_i \right )
\end{eqnarray*}
where $s_1, \ldots, s_{m+1}, s \in S$ are distinct, and $\xi'_s \in
M_n(s)$ is curvi-linear. Let $U_2$ be the closure of the subset of
$\hil{m+n, m+n+1} \times S \times \hil{m+1}$ consisting of all the
elements:
\begin{eqnarray}   \label{des_U2}
\left (\left (\sum_{i=1}^{m+1} s_i + \xi_s, \, \sum_{i=1}^{m+1} s_i
+ \xi'_s \right ), \,\, s, \,\, \sum_{i=1}^{m+1} s_i \right )
\end{eqnarray}
where $s_1, \ldots, s_{m+1}, s \in S$ are distinct, $\xi'_s \in
M_n(s)$ is curvi-linear, and $\xi_s \subset \xi'_s$. Note that
general  elements in $M_n(s)$ are curvi-linear. Moreover, since
$\xi'_s \in M_n(s)$ is curvi-linear, $\xi'_s$ uniquely determines
$\xi_s$ with $\xi_s \subset \xi'_s$.

A typical element $(\xi', s, \eta') \in Q^{[m+n+1,m+1]}$ is:
\begin{eqnarray*}
(s_1 + \ldots + s_{m+1} + \xi'_s, \,\, s, \,\, s_1 + \ldots +
s_{m+1})
\end{eqnarray*}
where $s_1, \ldots, s_{m+1}, s \in S$ are distinct, and $\xi'_s \in
M_n(s)$ is curvi-linear. Hence,
\begin{eqnarray}       \label{wa_comm_g-n.5}
h^*[Q^{[m+n+1,m+1]}] = [U_1] + n[U_2].
\end{eqnarray}
On the other hand, a typical element $\big ( (\xi,\xi'), s, (\eta,
\eta') \big ) \in \Wq^{[m+n,m]}$ is of the form:
\begin{eqnarray*}
\left (\left (\sum_{i=1}^m s_i + \xi'_s, \, \sum_{i=1}^{m+1} s_i +
\xi'_s \right ), \,\, s, \,\, \left (\sum_{i=1}^m s_i, \,
\sum_{i=1}^{m+1} s_i\right ) \right )
\end{eqnarray*}
where $s_1, \ldots, s_{m+1}, s \in S$ are distinct, and $\xi'_s \in
M_n(s)$ is curvi-linear. Hence,
\begin{eqnarray}       \label{wa_comm_g-n.6}
g_*[\Wq^{[m+n,m]}] = [U_1].
\end{eqnarray}
Combining (\ref{wa_comm_g-n.3}), (\ref{wa_comm_g-n.4}),
(\ref{wa_comm_g-n.5}) and (\ref{wa_comm_g-n.6}), we obtain
\begin{eqnarray}
& &g_{m+n+1}^*\mathfrak a_{-n}(\alpha)(A)
   - {\wa_{-n}(\alpha)}g_{m+1}^*(A)  \nonumber  \\
&=&n \cdot ({\w p}_1')_* \big ([U_2] \cdot
   (\rho \circ h)^*\alpha \cdot (p_2 \circ h)^*A \big ).
                                   \label{wa_comm_g-n.7}
\end{eqnarray}

Let $\phi = {\rm Id}_{\hil{m+n, m+n+1}} \times \rho_{m+1} \times
f_{m+1}$ be the morphism:
\begin{eqnarray*}
\hil{m+n, m+n+1} \times \hil{m+1, m+2} \to \hil{m+n, m+n+1} \times S
\times \hil{m+1}.
\end{eqnarray*}
From the descriptions (\ref{gen_pt_wq}) and (\ref{des_U2}), we
conclude that
\begin{eqnarray*}
[U_2] = \phi_*[\Wq_{m+1, n-1}].
\end{eqnarray*}
By (\ref{wa_comm_g-n.7}), the projection formula and
Lemma~\ref{wt^n}, we obtain
\begin{eqnarray*}
& &g_{m+n+1}^*\mathfrak a_{-n}(\alpha)(A)
   - {\wa_{-n}(\alpha)}g_{m+1}^*(A)     \\
&=&n \cdot ({\w p}_1' \circ \phi)_* \left ([\Wq_{m+1, n-1}]
   \cdot (\rho \circ h \circ \phi)^*\alpha \cdot
   (p_2 \circ h \circ \phi)^*A \right )    \\
&=&n \cdot \pi_{1*} \left ([\Wq_{m+1, n-1}]
   \cdot \pi_2^*\big ( \rho_{m+1}^*(\alpha)
   \cdot f_{m+1}^*(A) \big ) \right )    \\
&=&n \cdot \wt^{n-1} \big ( \rho_{m+1}^*(\alpha)
   \cdot f_{m+1}^*(A) \big )
\end{eqnarray*}
where $\pi_1$ and $\pi_2$ are the two projections of $\hil{m+n,
m+n+1} \times \hil{m+1, m+2}$.
\end{proof}
\section{\bf Incidence Hilbert schemes and Lie algebras}
\label{sect_loop}
\subsection{\bf The main Theorem}$\,$ With all the results obtained
in previous sections, we are ready to formulate and prove  the main
theorem.

\begin{definition}   \label{alg_h}
(i) Define ${\w {\mathfrak h}}_S$ to be the Heisenberg algebra
generated by the operators $\wa_n(\alpha)$, $n \in \Z$, $\alpha \in
H^*(S)$ and the identity operator ${\rm Id}_{\Wfock}$;

(ii) Define a Lie algebra structure on
\begin{eqnarray*}
\W{\mathfrak H}_S={\w {\mathfrak h}}_S \oplus H^*(S)\oplus\C\wt
\end{eqnarray*}
by declaring
%%change
\begin{eqnarray}   \label{alg_h.1}
[\mathfrak \wa(\alpha), \beta] = 0,\quad [\beta,\gamma]=0,\quad
[\wt, \mathfrak \wa(\alpha)]=0, \quad [\wt, \beta]=0
\end{eqnarray}
for operators $\mathfrak \wa(\alpha) \in {\w {\mathfrak h}}_S$ and
cohomology classes $\beta,\gamma \in H^*(S)$.
\end{definition}

%If $\mathfrak g=\mathfrak g_{-}\oplus \mathfrak g_0\oplus \mathfrak
%g_{+}$ is the triangular decomposition,  the lower triangular part
%of the loop algebra is given by the subalgebra (see \cite{Kac}):
%\begin{eqnarray*}
%\big(t^{-1}\C [ t^{-1}] \otimes_\C (\mathfrak g_0\oplus\mathfrak
%g_+)\big)\oplus \C[t^{-1}]\otimes_{\C }\mathfrak g_- .
%\end{eqnarray*}
%Let $\W{\mathcal L}_S$ be the lower triangular part of the loop
%algebra $\C[t, t^{-1}]\otimes_\C \W{\mathcal H}_S$.

\begin{theorem}  \label{thm_struc}
The space $\Wfock$ is a representation of the Lie algebra
$\W{\mathfrak H}_S$ with a highest weight vector being the vacuum
vector
\begin{eqnarray*}
\vac = 1_S \in H^0(S) = H^0(\hil{0,1})
\end{eqnarray*}
where  $1_S$ denotes the fundamental cohomology class of $S$.
\end{theorem}
\begin{proof}
By Lemma~\ref{wa_homo}, Proposition~\ref{prop_hei_inc},
Lemma~\ref{wt_homo}~(ii) and Proposition~\ref{comm_wtwa}, we see
that the space $\Wfock$ is a representation of the Lie algebra $
\W{\mathfrak H}_S$

Next, for each $0 \le i \le 4$, we fix a linear basis
\begin{eqnarray*}
\{ \alpha_{i, 1}, \ldots, \alpha_{i, b_i(S)} \}
\end{eqnarray*}
of $H^i(S)$. By Proposition~\ref{prop_hei_inc} and the existence of
the left inverse $\wt^\dagger$ of the operator $\wt$ in
Lemma~\ref{wt_homo}~(i), the following cohomology classes:
\begin{eqnarray}   \label{monomial}
& &\wt^m \cdot \prod_{i \text{ even}}
     \mathfrak a_{-n_{i, j_i}}(\alpha_{i, j_i})^{m_{i, j_i}}
   \cdot \prod_{i \text{ odd}}
     \mathfrak a_{-n_{i, j_i}}(\alpha_{i, j_i})^{m_{i, j_i}}
   \cdot \alpha_{k, j_k} \vac     \nonumber     \\
&=&\wt^m \cdot \prod_{i \text{ even}}
     \mathfrak a_{-n_{i, j_i}}(\alpha_{i, j_i})^{m_{i, j_i}}
   \cdot \prod_{i \text{ odd}}
     \mathfrak a_{-n_{i, j_i}}(\alpha_{i, j_i})^{m_{i, j_i}}
   \cdot \alpha_{k, j_k}
\end{eqnarray}
are linearly independent, where $0 \le i, k \le 4$, $0 \le j_i \le
b_i(S)$, $0 \le j_k \le b_k(S)$, $n_{i, j_i} > 0$, $m \ge 0$, $m_{i,
j_i} \ge 0$ for all the $i$ and $j_i$, and $m_{i, j_i} = 0$ or $1$
for odd $i$. Since the bi-degrees of the operators $\wt$ and
$\wa_{-n}(\alpha)$ are $(1, 2)$ and $(n, 2n-2+|\alpha|)$
respectively, we conclude from (\ref{Hnn1_bi}) that cohomology
classes in (\ref{monomial}) form a linear basis of $\Wfock$.
Therefore, the representation is the highest weight representation.
%%change
\end{proof}

\begin{example}
We express certain  distinguished  coholomogy classes by using
(\ref{monomial}).

First, consider the fundamental cohomology class of $\hil{n, n+1}$.
For $n \ge 0$, let
\begin{eqnarray}   \label{1_n}
{\bf 1}_{-n} = \frac{1}{n!} \, \mathfrak a_{-1}(1_S)^n, \qquad \w
{\bf 1}_{-n} = \frac{1}{n!} \, \wa_{-1}(1_S)^n.
\end{eqnarray}
For convenience, we also put ${\bf 1}_{-n} = 0$ and $\w {\bf 1}_{-n}
= 0$ when $n < 0$. Let $1_{\hil{n}} \in H^0(\hil{n})$ and
$1_{\hil{n, n+1}} \in H^0(\hil{n, n+1})$ be the fundamental
cohomology classes of $\hil{n}$ and $\hil{n, n+1}$ respectively.
Then it is well-known that
\begin{eqnarray}     \label{1_n.1}
1_{\hil{n}} = {\bf 1}_{-n} \vac
\end{eqnarray}
where by abusing notations, we also let $\vac \in H^0(\hil{0}) \cong
\C$ to stand for the cohomology class corresponding to the point
$\hil{0}$. By Lemma~\ref{wa_comm_f},
\begin{eqnarray}   \label{1_n.2}
1_{\hil{n, n+1}} = f_n^* 1_{\hil{n}} = f_n^* {\bf 1}_{-n} \vac = \w
{\bf 1}_{-n} f_0^*\vac = \w {\bf 1}_{-n} \vac.
\end{eqnarray}

Next, let $n \ge 1$ and $E_{n, n+1}$ be the exceptional divisor in
$\hil{n, n+1}$ with respect to the blowing-up morphism $\hil{n, n+1}
\to \hil{n} \times S$, i.e., $E_{n, n+1}$ is given by
\begin{eqnarray*}
E_{n, n+1} = \left \{ (\xi, \xi') \in \hil{n, n+1}|\, \Supp(\xi) =
\Supp(\xi') \right \}.
\end{eqnarray*}
From the definition of the operator $\wt$, we see that $[E_{1, 2}] =
\wt \vac$. Assume $n \ge 2$. Let $B_n$ be the boundary divisor in
the Hilbert scheme $\hil{n}$, i.e.,
\begin{eqnarray*}
B_n = \{ \xi \in \hil{n}|\, \#\Supp(\xi)\le n-1 \}.
\end{eqnarray*}
Then, $[B_n] = \mathfrak a_{-2}(1_S){\bf 1}_{-(n-2)} \vac$ and
$g_{n+1}^*[B_{n+1}] = f_n^*[B_n] + 2 [E_{n, n+1}]$. Thus,
\begin{eqnarray*}
   [E_{n, n+1}]
&=&\frac{1}{2} \, \left ( g_{n+1}^*[B_{n+1}] -
   f_n^*[B_n] \right )       \\
&=&\frac{1}{2} \, \left ( g_{n+1}^*
   \mathfrak a_{-2}(1_S){\bf 1}_{-(n-1)} \vac -
   \wa_{-2}(1_S) \w {\bf 1}_{-(n-2)} \vac \right ).
\end{eqnarray*}
From Proposition~\ref{wa_comm_g-n}~(ii), (\ref{1_n.1}) and
(\ref{1_n.2}), we conclude that
\begin{eqnarray}   \label{1_n.4}
   [E_{n, n+1}]
&=&\frac{1}{2} \, \left ( \wa_{-2}(1_S) g_{n-1}^*
   1_{\hil{n-1}} + 2 \wt f_{n-1}^*1_{\hil{n-1}} -
   \wa_{-2}(1_S) \w {\bf 1}_{-(n-2)} \vac \right )
   \nonumber  \\
&=&\wt \big ( 1_{\hil{n-1, n}} \big )  \nonumber  \\
&=&\wt \, \w {\bf 1}_{-(n-1)} \vac.
\end{eqnarray}
\end{example}

\subsection{\bf Quasi-projective case}
$\,$ We make some comments for a smooth quasi-projective surface
$S$. It has been proved in \cite{LQ1} that Cheah's formula
(\ref{Hnn1_bi}) holds for $S = \C^2$. However, it is unclear whether
(\ref{Hnn1_bi}) holds for an arbitrary quasi-projective surface $S$.
In the following, we will bypass this issue by assuming that
(\ref{Hnn1_bi}) holds for the quasi-projective surface $S$.

Since the surface $S$ is quasi-projective, some maps involved in the
definition of the Heisenberg operators, which are proper in the
projective case, are not proper anymore. Hence some modifications
are needed. We follow Nakajima's treatment of this issue in his book
\cite{Na2}. The main remedy is to choose appropriate (co)homology
theories so that the Heisenberg operators can still be defined.

Since $S$ is smooth, the Borel-Moore homology group $H^{lf}_i(S)$ is
dual to $H^{4-i}(S)$. It follows that the morphism (\ref{rho_n})
induces an $H^{lf}_i(S)$-module structure on:
\begin{eqnarray*}
\W{\mathbb H}_S=\bigoplus_{m=0}^{+\infty}H^*(\hil{m, m+1}).
\end{eqnarray*}

For the creation operators $\wa_{-n}(\alpha)$, where $n>0$, we take
$\alpha$ to be a class in the Borel-Moore homology group
$H^{lf}_*(S)$. For the annihilation operators $\wa_{n}(\beta)$,
where $n>0$, we take $\beta$ to be a class in the ordinary homology
group $H_*(S)$. Once we have this set-up, the Heisenberg operators
are well defined and the arguments in the previous sections are
still valid. Note that $H^{lf}_*(S)$ and $H_*(S)$ are dual to each
other. So we have the commutation relation as
(\ref{prop_hei_inc.1}).

The translation operator $\wt$ and its adjoint $\wt^\dagger$ can be
defined as usual without any change since the projections from the
subset $\Wq_m\subset \hil{m+1, m+2}\times \hil{m, m+1}$ to both
factors $\hil{m+1, m+2}$ and $\hil{m, m+1}$ are proper.

In summary, with this set-up, we have the same statement as in
Theorem~\ref{thm_struc} on $\W{\mathbb H}_S$. In \cite{LQ1}, we have
worked out this for $S = \C^2$ in the equivariant setting.

\end{document}